\documentclass[12pt]{article}

\usepackage{hyperref}
\usepackage{amssymb,amsmath,amsfonts,graphicx,subfigure,color,multicol,amsthm,float}
\usepackage[margin=1.0in]{geometry}
\newcommand{\Ss}{{\mathcal{S}}}

\newcommand{\Ff}{{\mathfrak{F}}}

\newcommand{\Cc}{{\mathbb{C}}}
\newcommand{\al}{\alpha}
\newcommand{\vn}{\mathbf{vN}}
\newcommand{\svn}{\mathbf{svN}}
\newcommand{\sch}{\mathbf{Sch}}
\newcommand{\lc}{\mathfrak{l c}\,}
\newcommand{\cc}{\mathfrak{c c}\,}
\newcommand{\Ee}{{\mathcal{E}}}

\newtheorem{theorem}{Theorem}[section]
\newtheorem{lemma}[theorem]{Lemma}

\newtheorem{remark}[theorem]{Remark}

\begin{document}

\title{Optimal subsets in the stability regions of multistep methods}



\author{Lajos L\'oczi\thanks{{\texttt{LLoczi@cs.elte.hu}}, Department of Numerical Analysis, E\"otv\"os Lor\'and University, and Department of Differential Equations, Budapest University of Technology and Economics, Hungary}}

\date{\today}



\maketitle

\begin{abstract}
In this work we study the stability regions of linear multistep or multiderivative
multistep methods for initial-value problems by using techniques that 
are straightforward to implement in modern computer algebra systems. In many applications, one is interested in (\textit{i}) checking whether a given subset of the complex plane (e.g.~a sector, disk, or parabola) is included in the stability region of the numerical method, (\textit{ii}) finding the largest subset of a certain shape contained in the stability region of a given method, or (\textit{iii}) finding the numerical method in a parametric family of multistep methods whose stability region contains the largest subset of a given shape.

First we describe a simple procedure to 
exactly calculate the stability angle $\alpha$ 
in the definition of $A(\alpha)$-stability by representing the root locus curve 
(RLC) of the multistep method as an implicit algebraic curve. 
As an illustration, we consider two finite families of implicit multistep methods.
We
exactly compute the stability angles for the $k$-step BDF methods
($3\le k\le 6$) and discover that the values of $\tan(\alpha)$ are 
surprisingly simple algebraic numbers of degree 2, 2, 4 and 2, respectively. 
In contrast, the corresponding values of $\tan(\alpha)$ for the $k$-step
second-derivative multistep methods of Enright ($3\le k\le 7$) are much more complicated;
the smallest algebraic degree here is 22.

Next we determine the exact value of the stability radius in the BDF family for each $3\le k\le 6$, that is, the radius of the largest disk in the left half of the complex plane, symmetric with respect to the real axis, touching the imaginary axis and lying in the stability region of the corresponding method. These radii turn out to be algebraic numbers of degree 2, 3, 5 and 5, respectively.

Finally, we demonstrate how some
Schur--Cohn-type theorems of recursive nature and \textit{not} relying on the RLC method can be used to  
exactly solve some optimization problems within infinite parametric families of multistep methods.
As an example, we choose a two-parameter family of implicit-explicit (IMEX) methods:
 we identify the unique method having the
largest stability angle in the family, then we 
find the unique method in the same family whose stability region contains
the largest parabola.
 

\end{abstract}

\section{Introduction}\label{introductionsection}


In the stability theory of one-step or multistep methods for initial-value problems, one is often interested in various geometric properties of the stability region $\Ss\subset \Cc$ of the method. 
In this work we study the shape of the stability region of  
linear multistep methods (LMMs) or multiderivative multistep methods (also known as
generalized LMMs) as follows. 

Suppose we are given 
 \begin{itemize}
\item[a)] a stability region $\Ss$, or
\item[b)] a family of stability regions $\Ss_\beta$ parametrized by some $\beta\in\mathbb{R}^d$,
\end{itemize}
and a family of subsets of $\Cc$, denoted by $\Ff$.  Due to their relevance in applications, we will consider the following three classes:
\begin{itemize}
\item $\Ff=\Ff^\text{\,sect}_\alpha$ is the family of infinite sectors in the left half of $\Cc$, with vertex at the origin, symmetric about the negative real axis, and parametrized by the sector angle $\alpha\in(0,\pi/2)$;
\item $\Ff=\Ff^\text{\,disk}_r$ is the family of disks in the left half of $\Cc$, symmetric with respect to the real axis, touching the imaginary axis, and parametrized by the disk radius $r>0$;
\item $\Ff=\Ff^\text{\,para}_m$ is the family of parabolas in the left half of $\Cc$, symmetric with respect to the real axis, touching the imaginary axis, and parametrized by some $m>0$.
\end{itemize}
Our goal is to find the set $H\in\Ff$ with the \emph{largest} parameter ($\alpha$, $r$, or $m$) such that \begin{itemize}
\item $H\subset\Ss$ in case a);
\item $H\subset \Ss_{\beta_\text{opt}}$ for some stability region in the family in case b), but $H\not\subset \Ss_{\beta}$ for $\beta\ne\beta_\text{opt}$.
\end{itemize}
We will present some tools to handle these shape optimization questions, and, as an illustration,  \textit{exactly} solve some of them by using \textit{Mathematica} version 11 in the BDF and Enright families (as LMMs and multiderivative multistep methods, respectively), and in an infinite family of IMEX methods with $d=2$ parameters.

\subsection{Motivation and main results}

\noindent  \fbox{\textbf{\,A\,}} When solving stiff ordinary differential equations, one desirable property of the numerical method is $A$-stability: a method is $A$-stable if the closed left half-plane $\{z\in\Cc : \mathrm{Re}(z)\le 0\}$
belongs to $\Ss$.
Many useful methods are not $A$-stable, 
still, $\Ss$ contains a sufficiently large infinite sector in the left
half-plane with vertex at the origin and symmetric
about the negative real axis. 
This leads to the notion of $A(\alpha)$-stability: a method is 
$A(\al)$-stable with some $0<\al<\pi/2$ if
\begin{equation}\label{Aalphadef}
\{z\in\Cc\setminus\{0\} : |\mathrm{arg}(-z)|\le \al\} \subset \Ss,
\end{equation}
where the argument of a non-zero complex number
satisfies $-\pi<\mathrm{arg}\le\pi$.
The largest $0<\al<\pi/2$ such that \eqref{Aalphadef} holds is referred to as the 
\textit{stability angle} 
of the method \cite{hairerwanner}. 
Various other stability concepts---such as $A(0)$-stability, $A_0$-stability, $\overset{\circ}{\text{A}}$-stability, stiff stability, or asymptotic 
$A(\al)$-stability---have also been defined, and theorems devised to test whether a given multistep method is stable in one of the above senses; see, for example,
\cite{liniger,cryer1973,jeltschmathcomp1976,jeltschnoteon,jeltschA0A0,jeltschmultidermultistep,jeltschcorr,friedlijeltsch,bickart,sanzserna}.  
There are various techniques to test $A(\alpha)$-stability for a given $\al$ value. 
In \cite{bickart}, for example, the sector on the left-hand side of \eqref{Aalphadef} is 
decomposed into an infinite union of disks, and a bijection between each disk and the left half-plane 
is established via fractional linear transformations to employ a Routh--Hurwitz-type criterion.
Another way of studying $A(\alpha)$-stability is to consider the root locus curve (RLC) of the multistep method \cite{hairerwanner}.
Based on the RLC and some theorems from complex analysis, \cite{norsett} 
presents a criterion for a LMM to be $A(\al)$-stable
for a given $\al$; the stability angle is then obtained as the solution of an optimization problem 
involving Chebyshev polynomials. 
The procedure in \cite{norsett} is formulated only for LMMs but not for multiderivative multistep methods. 

The first goal of the present work is to describe an elementary approach to 
exactly determine the stability angle of a LMM or multiderivative multistep method:
by eliminating the complex exponential function from the RLC and using a 
tangency condition, a system of polynomial equations in two variables is set up whose solution
yields the stability angle.
This process is easily implemented in computer algebra systems.
As an illustration, we 
consider two finite families: the BDF
methods \cite{gear,hairerwanner,norsett} as LMMs, and the second-derivative 
multistep methods of Enright \cite{enright,chakravartikamel,hairerwanner}.
With $\alpha_k^\text{BDF}$ denoting the stability angle of the $k$-step BDF method for $3\le k \le 6$, 
we show that $\tan\left(\alpha_k^\text{BDF}\right)$ 
is an unexpectedly simple algebraic number, having degree 2
for $k\in \{3, 4, 6\}$, and degree 4 for $k=5$;
see Table \ref{tab:1}. For the $k$-step Enright methods with $3\le k\le 7$,  
 the corresponding constants
$\tan\left(\alpha_k^\text{Enr}\right)$ (with approximate values listed in Table \ref{tab:2}) are much more complicated algebraic numbers
of increasing degree (starting with 22). 
As far as we know, exact values $\al\in (0,\pi/2)$ for the stability angles of multistep methods 
were not presented earlier in the literature.
\begin{remark}
The $k$-step BDF methods for $k\in \{ 1, 2\}$
are $A$-stable. For 
 $k\ge 7$ they are not zero-stable \cite{cryer,creedonmiller,hairerwannerarticle}, 
therefore not interesting from a practical point of view.
\end{remark}

\begin{remark}
In \cite[Table 1]{norsett} one finds some approximate values for the BDF stability angles, 
however, some of these values are  not correct.
The $k=3$ value is wrong because the polynomial $R_3$ is not computed properly.
The approximate values for $k=4$ and $k=5$ given in \cite[Table 1]{norsett} are correct (up to the given precision). The value for $k=6$ is again incorrect because an error
was committed in the minimization process.
If the optimization in \cite[Section 3]{norsett} is carried out
exactly with the correct $R_j$ polynomials, we recover the stability angle values  
in our Table \ref{tab:1}. The errors in \cite[Table 1]{norsett} propagated in 
the literature, see, for example, \cite[p. 242]{lambert}.   
As a consequence, some works that appeared 
in the current millennium also contain the erroneous angles.
In \cite[Chapter V.2, (2.7)]{hairerwanner} the correct approximate values are presented.
\end{remark}

\begin{remark}
At the time of writing this document, we learned (through personal communication) that \cite{akrivis} also contains the exact stability angles for the BDF methods with $3\le k \le 6$ steps: although they use a different technique to derive the results and the $\arcsin$ function to express the final constants,  the values  given in  \cite{akrivis} and our Table \ref{tab:1} are the same. Notice, however, that the stability angle for $k=5$ given in \cite{akrivis} has a slightly more complicated structure than the value in our Table \ref{tab:1}.
\end{remark}

\begin{remark}
The $k$-step Enright methods are $A$-stable again for $k\in \{ 1, 2\}$, see \cite{hairerwanner}, and unstable for $k\ge 8$.  More precisely, \cite{friedlijeltsch} proves that these methods
are not $A_0$-stable for $k\ge 8$, hence they cannot be stiffly stable either, see
\cite[Theorem 3]{jeltschmultidermultistep} (cf.~\cite{jeltschA0A0,jeltschcorr}).
However, in \cite[Chapter V.3, p. 276, Exercise 2]{hairerwanner} the stiff instability of the Enright formulas for $k\ge 8$ is still mentioned as an open problem.
\end{remark}

\noindent \fbox{\textbf{\,B\,}} The \textit{stability radius} of a multistep method is the largest number $r>0$ such that the inclusion
\[
\{z\in\mathbb{C} : |z+r|\le r\} \subset \Ss
\]
holds. The stability radius plays an important role when analyzing boundedness properties
of multistep methods. For example, it has been proved \cite[Theorem 3.1]{spijker2017} that this radius is the largest step-size coefficient for linear boundedness of a LMM satisfying some natural assumptions. 
\begin{remark} For LMMs (and for more general methods as well), various other step-size coefficients have been introduced in the context of linear
or non-linear problems. These coefficients govern the largest allowable step-size guaranteeing  
certain monotonicity or boundedness properties of the LMM, including the TVD and SSP properties \cite{David11}. These properties are relevant, for example, in the time integration of method-of-lines semi-discretizations of hyperbolic conservation laws \cite{spijker2013,hundsdorferspijkermozartova2012,LLexact}.
 \end{remark}
 \begin{remark}
 In \cite{DKKATLL} the largest inscribed and smallest circumscribed (semi)disks are computed for certain one-step methods.
\end{remark}

The second goal of the present work is to compute the {stability radius} for some multistep methods. We will achieve this by using again the algebraic form of the RLCs.  
Table \ref{tab:3} contains the exact values in the BDF family for $3\le k \le 6$.\\

\noindent \fbox{\textbf{\,C\,}} The RLC, as the graph of a $[0,2\pi]\to \mathbb{C}$ function (or a union of such functions for generalized LMMs), yields information about the 
boundary of the stability region, $\partial \Ss$. It is known, however, that 
in general the RLC does not coincide with $\partial \Ss$
(see 
Figure \ref{trueboundary}). 
This does not pose a problem when 
a fixed multistep method is considered---one can 
evaluate the roots of the characteristic polynomial at
finitely many test points sampled from different components of $\Cc$ determined by the RLC
to see which component belongs to $\Ss$ and which one to $\Cc\setminus \Ss$. 
But when working with parametric 
families of multistep methods, the precise identification of the stability region boundaries
or components 
can become challenging with the RLC method.
One can overcome this difficulty for example by invoking a reduction process, the Schur--Cohn reduction, 
formulated in e.g.~\cite{miller}.
Instead of using auxiliary 
fractional linear transformations and applying Routh--Hurwitz-type criteria
\cite{morrismarden,jeltschproc} as mentioned above,
these Schur--Cohn-type theorems in \cite{miller} are directly tailored to the 
context of multistep methods to locate the 
roots of the characteristic polynomials with respect to the unit disk. 

The third goal of the present work is to demonstrate the 
effectiveness of the Schur--Cohn reduction when we
solve two optimization case studies in
a family of implicit-explicit (IMEX) multistep methods taken from \cite{hundsdorferruuth}.
On the one hand, we find the method in the IMEX family that has the largest stability angle,
that is, the method whose stability region contains the largest sector (see our Theorem \ref{wedgethm}). 
On the other hand, 
we illustrate the versatility of the reduction technique by also finding the method
whose stability region contains the largest parabola (see Theorem \ref{parabolathm}); the inclusion of a parabola-shaped region in $\Ss$ is relevant when studying semi-discretizations of certain partial differential equations (PDEs) of  advection-reaction-diffusion type
\cite{jeltschproc,calvo,hundsdorferruuth}. The chosen IMEX family is described by two real parameters, and the corresponding characteristic polynomial is cubic. 
 The Schur--Cohn reduction process recursively decreases    
the degree of the characteristic polynomial, so instead of analyzing the \textit{roots} of high-degree polynomials, we finally need to check polynomial inequalities in the parameters present in the 
\emph{coefficients} of the original polynomial.  
Besides the two real parameters, two complex variables are involved 
in our calculations---the non-trivial interplay between these six real variables determines
 the optimum in both cases.
We emphasize that we solve the optimization problem exactly, and RLCs are
not relied on in the rigorous part of the proofs (only when setting up conjectures about the optimal values). 
\begin{remark}
The Schur--Cohn reduction is also used in \cite{jeltschbrown} to explore certain properties of a discrete parametric
family of multistep methods.  Conditions for disk or segment inclusions in the stability regions of a two-parameter 
family of multistep methods are formulated in \cite{oliveira}.
Optimality questions about the size and shape of the stability regions of one-step or multistep methods are investigated in detail in \cite{jeltschnevanlinna}. Properties
of optimal stability polynomials and stability region optimization in parametric families of one-step methods are discussed, for example, in \cite{jeltschtorrilhon,ketchesonahmadia}. 
\end{remark}

\subsection{Structure of the paper}

In Section \ref{prelimsection}, we introduce some notation. 
In Sctions \ref{schurcohnsection}--\ref{sectionvanishingcoeffs}, we review the Schur--Cohn reduction and the definition of the stability region
of a multistep method.
In Sections \ref{RLCLMMsection}--\ref{RLCMMMsection}, the
definition of the root locus curve is recalled in two special cases: for linear multistep methods and 
for second-derivative multistep methods.
Here we consider the BDF and Enright families as concrete examples. 

Regarding the new results, a simple algebraic technique is described in Section \ref{section23} to exactly compute the stability angle of a linear multistep or multiderivative multistep method. Stability angles for the BDF and Enright families are tabulated in Sections \ref{BDFsection}--\ref{Enrightsection}.
In Section \ref{diskinclusionsection}, we exactly compute the stability radii in the BDF family by using
the same approach. In Section \ref{section7}, we first describe a 
two-parameter family of IMEX multistep methods, in which we determine the unique method with
the largest stability angle, then, in Section \ref{section8}, the unique method 
whose stability region contains the largest parabola. The techniques
in Sections \ref{section7}--\ref{section8} do \textit{not} rely on root locus curves but use the Schur--Cohn reduction instead; the full proofs are deferred to Appendices \ref{appendixA} and \ref{appendixB}. 

\section{Preliminaries}

\subsection{Notation}\label{prelimsection}

The set of natural numbers $\{0,1,\ldots\}$ is denoted by $\mathbb{N}$. 
For $z\in\mathbb{C}$, $\mathrm{Re}(z)$, $\mathrm{Im}(z)$, and
$\overline{z}$  denote the real and imaginary parts, and the
conjugate of $z$, respectively, and $i$ is the imaginary unit. 
The boundary of a (possibly unbounded) set $H\subset\mathbb{C}$
is $\partial H\subset\mathbb{C}$. 
When describing certain algebraic numbers of higher degree, a polynomial $\sum_{j=0}^n a_j x^j$ 
with $a_j\in\mathbb{Z}$, $a_n\ne 0$ and $n\ge 3$ will be represented simply by its coefficient list $\{a_n, a_{n-1},\ldots, a_0\}$. For a polynomial $Q(z)=\sum_{j=0}^{n}a_j z^j$ with $0\le n
\in \mathbb{N}$, 
$a_j\in\mathbb{C}$ ($0\le j\le n$) and $a_n\ne 0$, we denote its
degree, leading coefficient and constant coefficient by $\deg Q=n$,    
 $\lc Q=a_n$ and $\cc Q=a_0$. The acronyms RLC and LMM stand for {root locus curve} and {linear multistep method}, respectively.

 \subsection{The Schur--Cohn reduction}\label{schurcohnsection}

In the rest of this section we assume that $Q$ is a univariate polynomial with $\deg Q\ge 1$, 
and follow the terminology of \cite{miller}---we have explicitly added the $\deg Q\ge 1$ condition,  
being implicit in \cite{miller}. We say that
 \begin{itemize}
\item $Q$ is a \textit{Schur polynomial}, $Q\in\sch$, if
its roots lie in the open unit disk;
\item $Q$ is a \textit{von Neumann polynomial}, $Q\in\vn$, if
its roots lie in the closed unit disk;
\item $Q$ is a \textit{simple von Neumann polynomial}, $Q\in\svn$, if
$Q\in\vn$ and roots with modulus 1 are simple.
\end{itemize}
\begin{remark}
The class $\sch$ is referred to as \emph{strongly stable polynomials} in
 \cite[p. 345]{butcher}.
\end{remark}
\begin{remark}
The property  $Q\in\svn$ is often expressed by saying that  $Q$ satisfies the \emph{root condition}.
\end{remark}
The \textit{reduced polynomial} of $Q(z)=\sum_{j=0}^{n}a_j z^j$ is defined as
\[
Q^{\mathbf{r}}(z):=\frac{\overline{a_n}\cdot \left(\sum_{j=0}^{n}a_j z^j\right)-a_0\cdot\left(\sum_{j=0}^{n}\overline{a_{n-j}}z^j\right)}{z}=
\]
\[
\sum_{j=1}^{n}\left(\overline{a_n} \cdot a_j
-a_0\cdot \overline{a_{n-j}}\right) z^{j-1},
\]
so we have $\deg Q^{\mathbf{r}}\le (\deg Q)-1$.
When this reduction process is iterated, we write $Q^{\mathbf{r}\mathbf{r}}$ for $\left(Q^\mathbf{r}\right)^\mathbf{r}$, for example. The following theorems 
from \cite{miller} use the notion of the reduced polynomial and the derivative
to formulate necessary and sufficient conditions for a polynomial to
be in the above classes. In all three theorems below it is assumed that $\lc Q \ne 0 \ne \cc Q$ and
$\deg Q\ge 2$. 
\begin{theorem}\label{schthm} $Q\in\sch \Leftrightarrow (|\lc Q|>|\cc Q| \text{ and } Q^\mathbf{r}\in\sch)$.
\end{theorem}
\begin{theorem}\label{vnthm} $Q\in\vn \Leftrightarrow \text{ either } (|\lc Q|>|\cc Q| \text{ and } Q^\mathbf{r}\in\vn) \text{ or } (Q^\mathbf{r}\equiv 0 \text{ and } Q'\in\vn)$.
\end{theorem}
\begin{theorem}\label{svnthm} $Q\in\svn \Leftrightarrow \text{ either } (|\lc Q|>|\cc Q| \text{ and } Q^\mathbf{r}\in\svn) \text{ or } (Q^\mathbf{r}\equiv 0 \text{ and } Q'\in\sch)$. 
\end{theorem}

\begin{remark}
Let us consider the following example when applying the theorems above, e.g.~Theorem 
\ref{vnthm}. 
For any $\lambda>0$, we set
$Q_\lambda(z):=z ^2+\lambda  i z +1$. 
Then the roots of $Q_\lambda$ satisfy $|z_1(\lambda)|<1<|z_2(\lambda)|$, so $Q_\lambda\notin \vn$, and  
$Q_\lambda^\mathbf{r}=2\lambda i$. 
This shows that it can happen that
the degree of the original polynomial is $>1$, but its reduced polynomial is a non-zero
constant, so the relation $Q^\mathbf{r}\in\vn$ is \emph{undefined}. In these cases, when $Q^\mathbf{r}$ is a 
non-zero constant, notice that neither  $|Q^\mathbf{r}|<1$, nor $|Q^\mathbf{r}|=1$, nor 
$|Q^\mathbf{r}|>1$ can help us in general to 
 determine whether $Q\in \vn$ or not (of course, the other
condition $|\lc Q|>|\cc Q|$ is violated now); cf.~the sentence above \cite[Theorem 5.1]{miller}. 
\end{remark}

\subsection{The stability region of a multistep method}\label{sectionvanishingcoeffs}


Stability properties of a broad class of numerical methods (including 
Runge--Kutta methods, linear multistep methods, or multiderivative multistep methods)
for solving initial value problems of the form 
\begin{equation}\label{IVP}
y'(t)=f(t,y(t)), \quad y(t_0)=y_0
\end{equation} 
can be
analyzed by studying the stability region of the method. When an $s$-stage $k$-step 
method ($s\ge 1$, $k\ge 1$ fixed positive integers; for $k=1$ we have a one-step method, while for $k\ge 2$ a multistep method) with constant
step size $h>0$ is applied to
the linear test equation $y'=\lambda y$ ($\lambda\in \Cc$ fixed, $y(0)=y_0$ given), the method yields a numerical solution $(y_n)_{n\in\mathbb{N}}$ that
approximates the exact solution $y$ at time $t_n:=t_0+n h$
and
satisfies a recurrence relation of the form \cite{jeltschnevanlinna}
\begin{equation}\label{recurrence}
\left\{
\begin{aligned}
   & \sum_{j=0}^s \sum_{\ell=0}^k a_{j, \ell} \, \mu^j \, y_{n+\ell} =0, \quad \ n\in\mathbb{N}, \\ 
   &   a_{j, \ell}  \in \mathbb{R}, \ \  \sum_{j=0}^s |a_{j,k}|>0, \ \  \mu:=h\lambda.
\end{aligned}
\right.
\end{equation}
 The characteristic polynomial associated with the method takes the form
\begin{equation}\label{phidef}
\Phi(\zeta,\mu):=   \sum_{j=0}^s \sum_{\ell=0}^k a_{j, \ell} \, \mu^j \, \zeta^\ell \quad
(\zeta\in\Cc).
\end{equation}
With $\Phi(\cdot,\mu)$ abbreviating the polynomial $\zeta\mapsto\Phi(\zeta,\mu)$, the stability region of the method is defined as
\begin{equation}\label{stabregdef}
\Ss:=\{ \mu\in\mathbb{C} : \text{the degree of } \Phi(\cdot,\mu)\text{ is exactly } k, \text{ and } \Phi(\cdot,\mu)\in\svn\}.
\end{equation}






\begin{remark}\label{nonvanishingremark}
Some other variations of the above definition of the stability region of a multistep method have also been proposed in the literature, see, e.g.~\cite{jeltschmulti}. In \cite[p. 344]{butcher}, the ``open stability region'' is defined as the set
\[
\{ \mu\in\mathbb{C} :  \Phi(\cdot,\mu)\in\sch\},
\]
see also \cite[p. 348]{suli}, \cite[p. 452]{gautschi} or \cite{lambert}. 
In e.g.~\cite{hairerwanner,kirlinger}, the stability region of the method \eqref{recurrence} is defined as
\begin{gather}\label{hairerstabilityregion}
\{ \mu\in\mathbb{C} : \text{all roots } \zeta_j(\mu) \text{ of } \zeta \mapsto \Phi(\zeta,\mu)
 \text{ satisfy } |\zeta_j(\mu)|\le 1, \\
\text{ and multiple roots satisfy } |\zeta_j(\mu)|<1\},\nonumber
\end{gather}
that is, essentially, $ \Phi(\cdot,\mu)\in\svn$.
In \cite[Formula (2.5)]{jeltschnevanlinna} the stability region is given by
\begin{gather}\label{extendedstabilityregion}
\{ \mu\in\overline{\mathbb{C}} : \text{roots } \zeta_j \text{ of \,} \Phi(\zeta,\mu)=0
 \text{ satisfy } |\zeta_j(\mu)|\le 1, \\
\text{ and if } |\zeta_j|=1, \text{\,then it is a simple root}\},\nonumber
\end{gather}
with $\overline{\mathbb{C}}$ denoting the extended complex plane. 

We can regroup the terms in \eqref{phidef} as $\Phi(\zeta,\mu)=\sum_{\ell=0}^k C_\ell(\mu)\zeta^\ell$
with some suitable polynomials $C_\ell$. The inequality condition in \eqref{recurrence} implies that the leading coefficient $C_k$ does not vanish \emph{identically}; it may happen that for some exceptional $\mu$ values the leading coefficient is zero:
\[
\Ee:=\{\mu\in\Cc : C_k(\mu)=0\}.
\]
For example, for the implicit Euler (IE) method $\Phi(\zeta,\mu)=\varrho(\zeta)-\mu\sigma(\zeta)=(1-\mu)\zeta-1$ with $\varrho(\zeta):=\zeta-1$ and $\sigma(\zeta):=\zeta$, so $\Ee=\{1\}$. For the $2$-step BDF method (BDF2), 
$\Phi(\zeta,\mu)=(3-2\mu)\zeta ^2-4 \zeta +1$, hence $\Ee=\{3/2\}$. If definition 
\eqref{hairerstabilityregion} (or \eqref{extendedstabilityregion}) is interpreted \emph{formally}, we have for the IE method that  
$\Ee=\{1\}\subset\Ss$ (because \eqref{hairerstabilityregion} is satisfied vacuously). 
Similarly, for the BDF2 method, $\Ee=\{3/2\}\subset\Ss$ (because then the 
unique root of $\Phi(\zeta,3/2)=0$ is $\zeta=1/4$).

However, elements of $\Ee$ or $\Ee\cap \Ss$ can be problematic.\\
(i) For $\mu\in\Ee$,  the order of the recursion \eqref{recurrence} decreases, thus, in general, the starting values $y_0, y_1, \ldots, y_{k-1}$ of the numerical method cannot be chosen arbitrarily.\\ 
(ii) Some exceptional values $\mu\in\Ee\cap \Ss$ can be located in the interior of the corresponding region of \emph{instability} of the method---this is the case for example for both the IE and BDF2 methods. 
When the step size $h > 0$ is chosen in a way that $\mu\in\Ee\cap \Ss$ is such an isolated value,
the recursion \eqref{recurrence} generated by the numerical method becomes practically useless (it
quickly ``blows up'' for arbitrarily small perturbations of $h$).\\
(iii) RLCs are often
used to identify the boundary $\partial \Ss$ of the stability region (see Sections \ref{RLCLMMsection}--\ref{RLCMMMsection} below). In \cite[Definition (2.21)]{jeltschnevanlinna}, the RLC is given by 
 \[\Gamma:=\{ \mu\in \overline{\Cc} : \exists \zeta\in\Cc\text{ with } |\zeta|=1 \text{ and } 
\Phi(\zeta,\mu)=0\}.\] 
It can happen that $\partial \Ss$ is a proper subset of the corresponding RLC (see, for example, our Figure \ref{trueboundary}), but in \cite[Corollary 2.6]{jeltschnevanlinna} it is shown that for a numerical method satisfying \emph{Property C} (see \cite[Formula (2.9)]{jeltschnevanlinna} or \cite[Definition 4.7]{hairerwanner}), the RLC \emph{coincides} with $\partial \Ss$. 
According to \cite[Section V.4]{hairerwanner}, all one-step methods have Property C, so the IE method
also has. 
And indeed, applying \cite[Proposition 2.7]{jeltschnevanlinna} to the IE method we now have that $\varrho$ and $\sigma$ have no common root and $\varrho/\sigma$ is univalent on the set 
$\{ z\in\overline{\mathbb{C}} :  |z-1|>1\}$, so $Q(\mu)=1/(1-\mu)$ has Property C.
Thus for the IE method $\partial \Ss=\Gamma$. As we have seen above, $1\in\Ee\cap\Ss$, so
$1\in \partial\Ss$. On the other hand, $\Phi(\zeta,1)=\varrho(\zeta)-\sigma(\zeta)=-1$, so
$1\notin\Gamma=\partial\Ss$. This apparent contradiction seems to indicate that the authors of \cite{jeltschnevanlinna} interpreted 
definition \eqref{extendedstabilityregion} \emph{intuitively}: a root $\zeta=\infty$ is tacitly introduced
as soon as the leading coefficient $C_k(\mu)$ becomes zero. So \cite[Corollary 2.6]{jeltschnevanlinna},
for example, actually relies on definition \eqref{stabregdef} rather than on definition \eqref{extendedstabilityregion} (or
\eqref{hairerstabilityregion}).

The problem of vanishing leading coefficient is implicitly avoided in \cite[p.~66]{lambert}, or in \cite{sanzserna}, because they impose a requirement on ``all the roots $r_s$ ($s=1, \ldots, k$)''.
Definition \eqref{stabregdef} above with the non-vanishing leading coefficient
essentially appears, for example, in \cite[Section 2.1]{spijker2013} 
(where it is formulated for LMMs, that is, for $s=1$ in \eqref{recurrence}), or in \cite[Section 2]{spijker2017}. 

Notice that, with the theorems cited in our Section \ref{schurcohnsection}, one can directly
investigate the stability region of a numerical method, \emph{without} constructing the corresponding RLC 
or without analyzing the relation between $\partial \Ss$ and the RLC (see Sections \ref{section7}--\ref{section8} below). 

Finally we remark that the above considerations also play an important role, e.g.~in control theory \cite[Chapter 1]{bhatta}, where a ``degree invariance" (i.e., ``no degree loss'') condition is incorporated in the \emph{Boundary Crossing Theorem}. \cite[Chapter 1]{bhatta} also recalls several stability results for polynomials, e.g.~the Routh--Hurwitz,  Jury, or  the recursive Schur(--Cohn) stability tests.
\end{remark}


\subsection{The RLC of a LMM}\label{RLCLMMsection}

A linear multistep method for \eqref{IVP}
has the form
\begin{equation}\label{multistepdefn}
\sum_{j=0}^k (\alpha_j y_{n+j}-h \beta_j f_{n+j})=0,
\end{equation}
where 
$f_m:=f(t_m, y_m)$, 
and the numbers $\al_j\in\mathbb{R}$ and $\beta_j\in \mathbb{R}$ ($j=0, \ldots, k$)
 are the suitably chosen method coefficients with $\alpha_k\ne 0$. 
The method is implicit, if $\beta_k\ne 0$.
By setting
\[\varrho(\zeta):=\sum_{j=0}^k \al_j \zeta^j \quad \text{and}\quad \sigma(\zeta):=\sum_{j=0}^k \beta_j \zeta^j,\]
the associated characteristic polynomial \eqref{phidef} becomes
\begin{equation}\label{P1def}
\Phi(\zeta,\mu)\equiv P_1(\zeta,\mu):=\varrho(\zeta)-\mu \sigma(\zeta).
\end{equation}
One way to study the stability region \eqref{stabregdef}, or its boundary $\partial \Ss$ in the complex plane is to depict the RLC corresponding to the method \cite{hairerwanner}: observe that $P_1$ is linear in $\mu$, so $P_1(\zeta,\mu)=0$ implies $\mu=\varrho(\zeta)/\sigma(\zeta)$ (for $\sigma(\zeta)\ne 0)$. The
RLC is then the image of the parametric curve 
\begin{equation}\label{muasaratio}
[0,2\pi]\ni \vartheta\mapsto \mu(\vartheta):=\frac{\varrho\left(e^{i\vartheta}\right)}{\sigma\left(e^{i\vartheta}\right)}.
\end{equation}

\subsubsection{RLCs for the BDF methods}\label{section21}

\begin{figure}[H]
\centerline{\includegraphics[width=.58\textwidth]{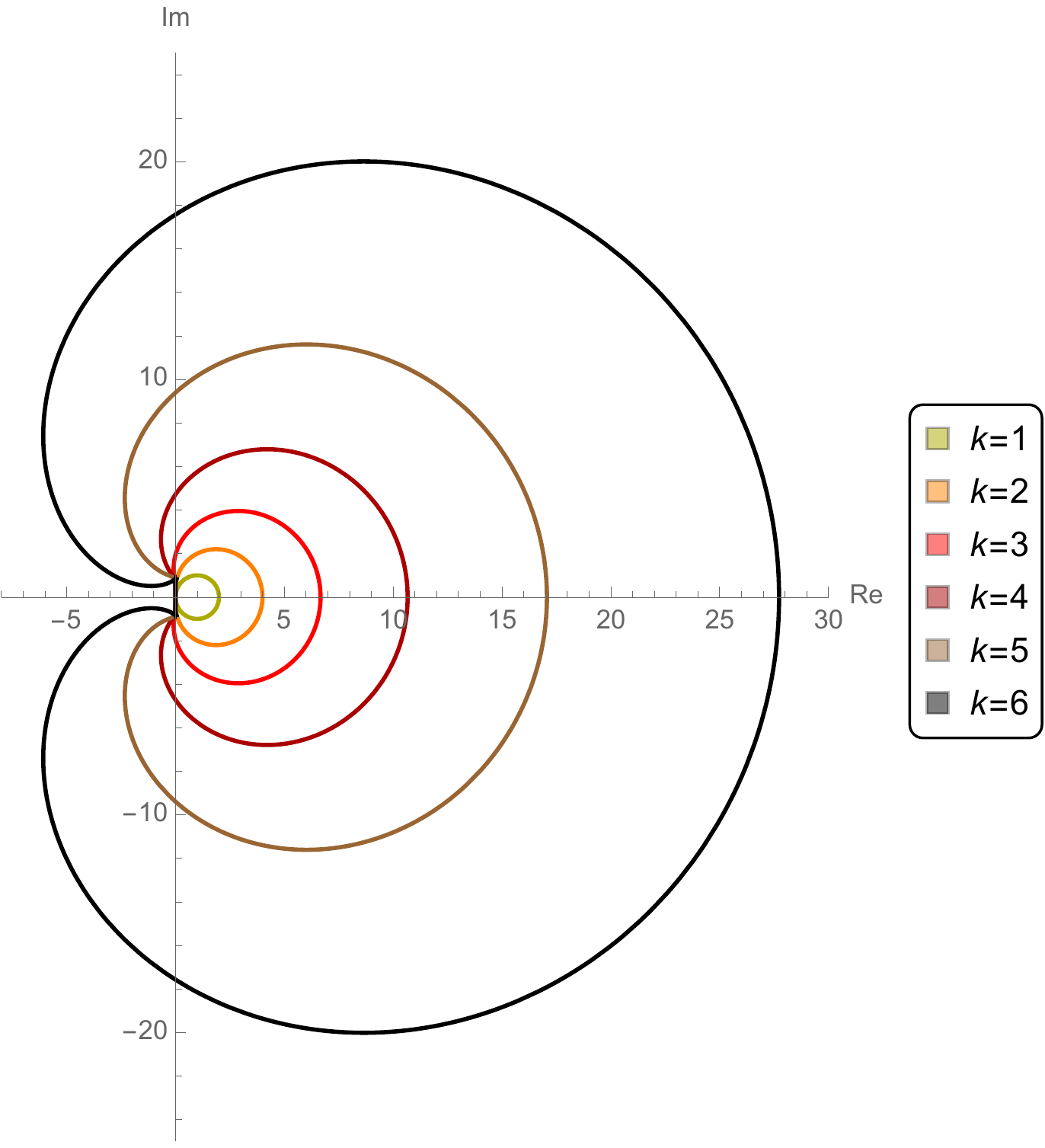}}
\caption{RLCs for the $k$-step BDF methods for $1\le k\le 6$. The stability
region of the method in each case is the unbounded component of $\Cc$.}
\label{fig:1}      
\end{figure}

\begin{figure}
\subfigure{
\includegraphics[width=0.435\textwidth]{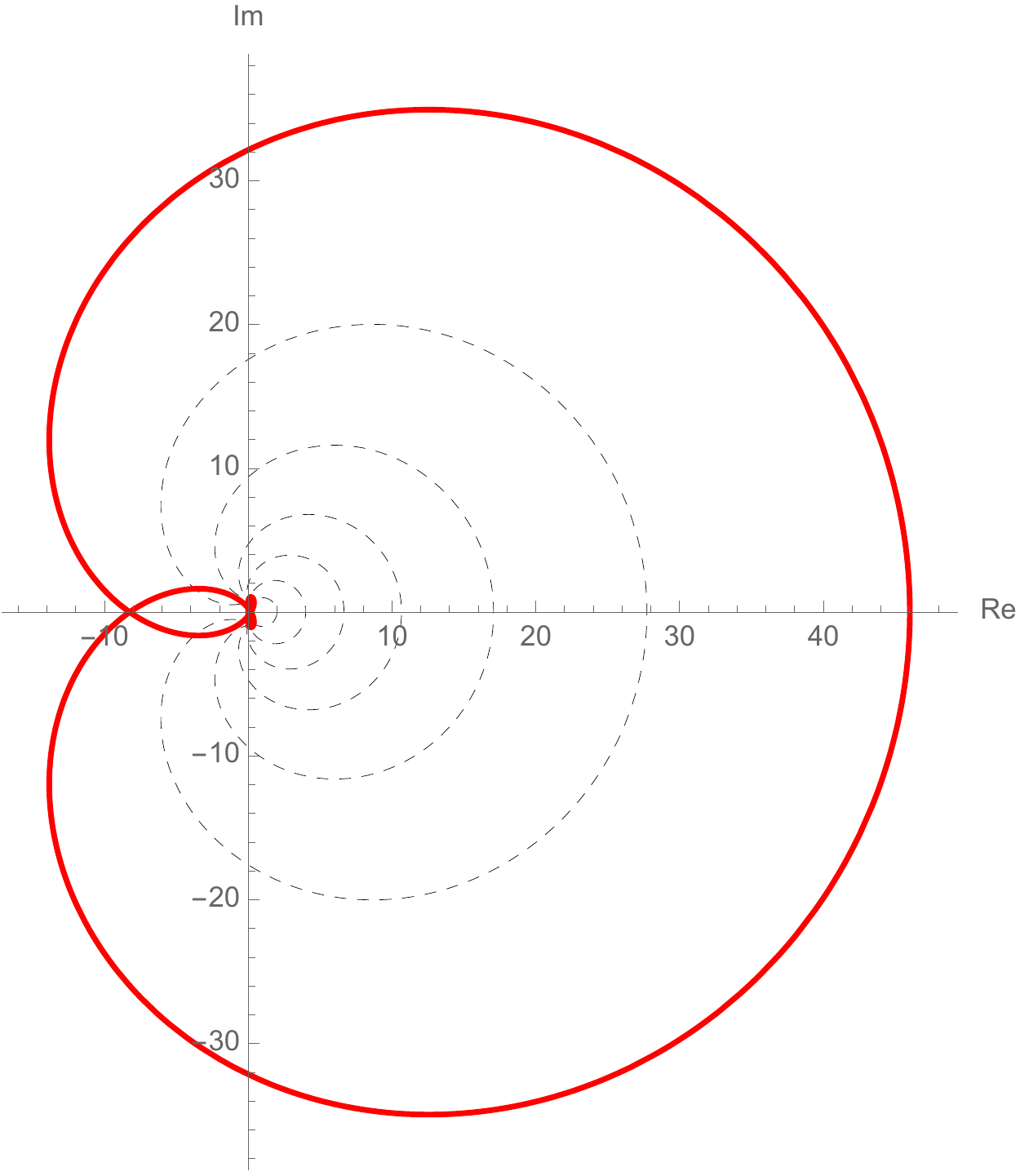}}
\subfigure{
\includegraphics[width=0.56\textwidth]{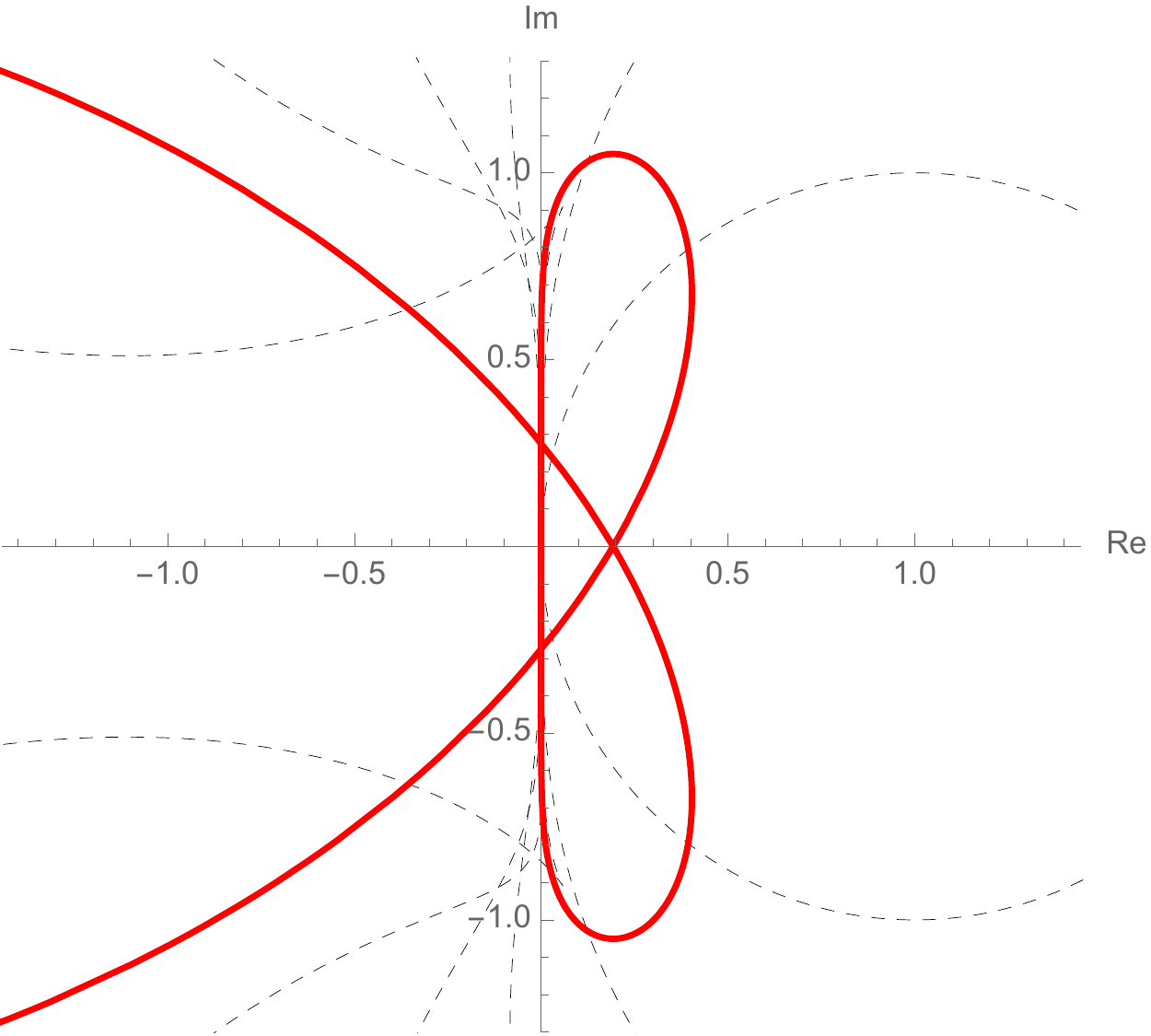}}
\caption{RLC for the unstable $7$-step BDF method in red (left), and a close-up near the origin (right). For comparison, 
the curves from Figure \ref{fig:1} are also superimposed as dashed gray curves. 
\label{splitfigure}}
\end{figure}

\begin{figure}
\subfigure{
\includegraphics[width=0.48\textwidth]{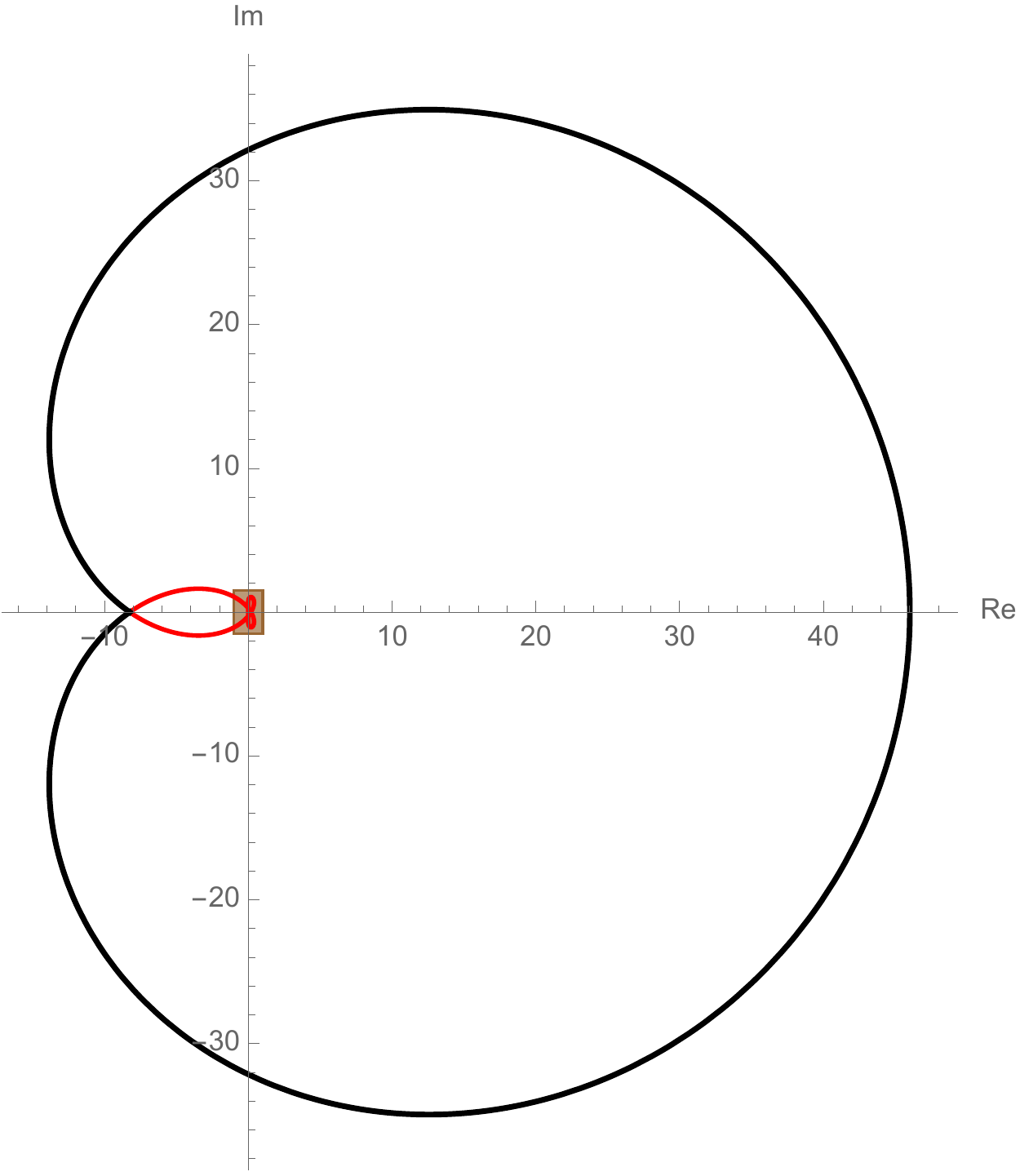}}
\subfigure{
\includegraphics[width=0.5\textwidth]{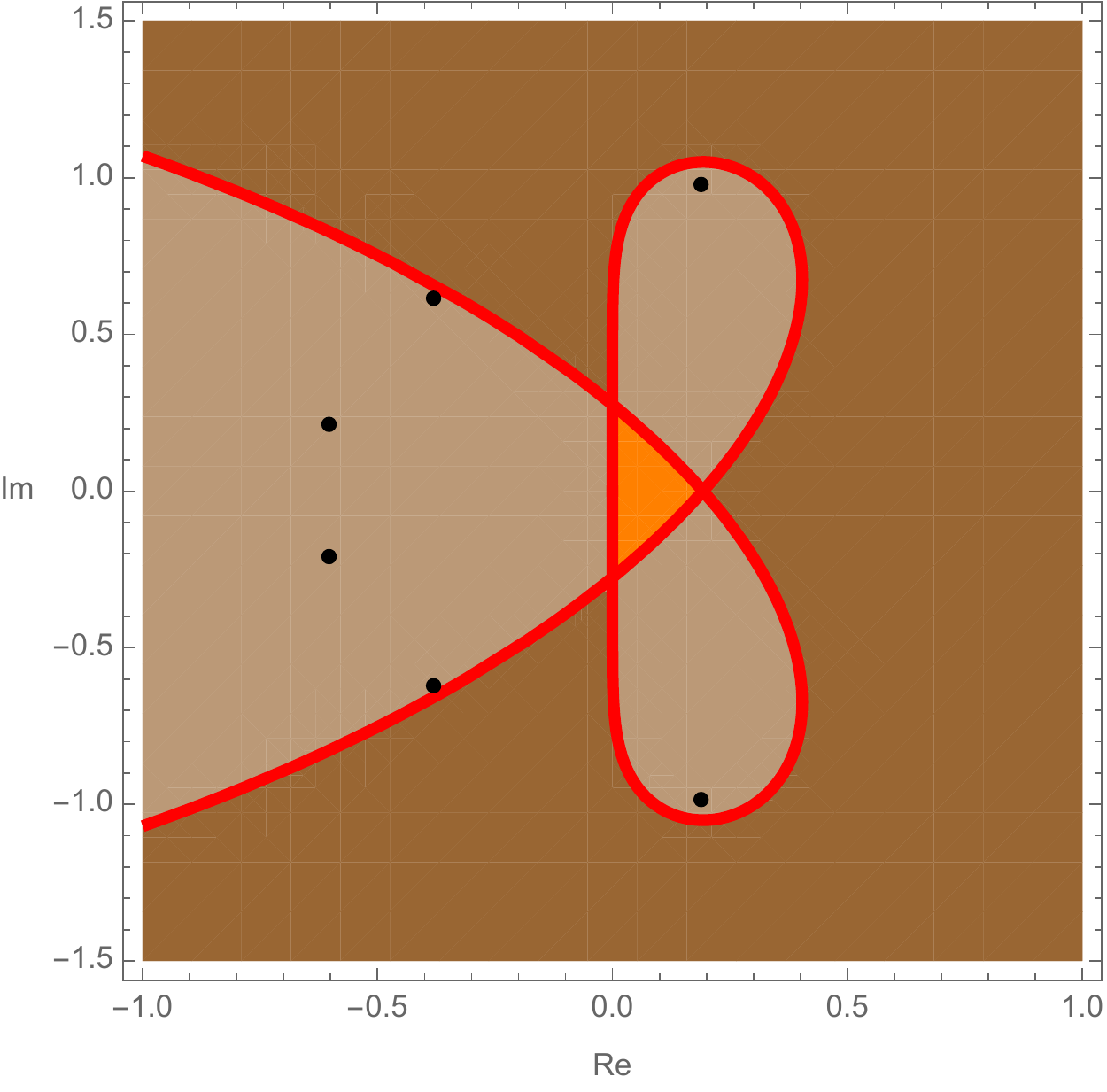}}
\caption{The black curve in the left figure shows the boundary $\partial \Ss$ of the stability region  of the (unstable) $7$-step BDF method; $\partial\Ss$ is non-differentiable at one point. The stability region is the unbounded outer component. The red curve segment
near the origin is not part of $\partial \Ss$, it is a subset only of the RLC as
displayed in Figure \ref{splitfigure}. The small brown rectangle in the center is shown
in detail in the right figure. The red curve in the right figure is again the RLC.
The 6 black dots depict the set of $\mu$ values
such that $P_1(\cdot,\mu)$ in \eqref{P1def} has multiple roots
(there are no other $\mu\in\mathbb{C}$ parameters with this property for $k=7$). The polynomial $P_1(\cdot,\mu)$
has 1, 2 and 3 roots outside the unit disk for $\mu$ values in the dark brown, 
light brown and orange regions, respectively; $P_1(\cdot,\mu)$ cannot have 4 or more roots
outside the unit disk. Each of the three self-intersections of the
RLC in this figure (as well as the self-intersection of the RLC 
 seen only in the left figure) 
corresponds to a $\mu$ value for which $P_1(\cdot,\mu)$
has two distinct roots with modulus 1.
Exactly computing, for example, the unique value of 
$\mu_\dagger\approx -2.68886\cdot 10^{-6}+ 0.275988 i$ in the open upper half-plane  where the RLC crosses itself was a non-trivial
task: it took \textit{Mathematica} 86 minutes to explicitly determine 
the coefficients of the integer polynomial defining $\mu_\dagger$ and having degree 30.
The RLCs for the $k$-step BDF methods with $1\le k \le 6$ do not have any self-intersections;
 other singularities may occur, see Figure \ref{BDF_6_cusp}.
 \label{trueboundary}}
\end{figure}

Each member of the BDF family is a special case of \eqref{multistepdefn}. The $k$-step BDF method
(having order $k$) is given by
\[
\sum_{j=1}^k \frac{1}{j}\nabla^j y_{n+1}=h f_{n+1},
\]
where $\nabla$ denotes the backward difference operator $\nabla y_{n+1}:=y_{n+1}-y_n$, and
$\nabla^j y_{n+1}:=$ $\nabla^{j-1} y_{n+1}-$ $\nabla^{j-1} y_n$ (for $j>1$). It is known
\cite{hairerwanner} that the corresponding RLC is
\begin{equation}\label{BDFconcreteformula}
\mu(\vartheta)\equiv \sum_{j=1}^k \frac{1}{j}(1-e^{-i\vartheta})^j.
\end{equation}
Figures \ref{fig:1}--\ref{trueboundary} 
show the RLCs for some BDF methods.

\subsection{The RLC of a multiderivative multistep method}\label{RLCMMMsection}

A second-derivative multistep method 
is more general than \eqref{multistepdefn} and can be written as
\begin{equation}\label{2nddermultistepdefn}
\sum_{j=0}^k (\alpha_j y_{n+j}-h \beta_j f_{n+j}-h^2 \gamma_j g_{n+j})=0,
\end{equation}
where $g_n:=g(t_n,y_n)$ with $g(t,y):=\partial_1 f(t,y)+\partial_2 f(t,y)\cdot f(t,y)$, and 
the method is determined by the coefficients $\al_j$, $\beta_j$ and $\gamma_j$, see  \cite{hairerwanner}.
Now the associated characteristic polynomial \eqref{phidef} becomes
\[
\Phi(\zeta,\mu)\equiv P_2(\zeta,\mu):=\sum_{j=0}^k (\al_j -\mu \beta_j-\mu^2 \gamma_j)\zeta^j.
\]
This time we have two RLCs: 
\begin{equation}\label{mu12asaquadratic}
[0,2\pi]\ni\vartheta\mapsto \mu_{1,2}(\vartheta),
\end{equation}
where $\mu_{1,2}$ are the two solutions of $P_2\left(e^{i\vartheta},\mu\right)=0$.
For any choice of the method coefficients $\al_j$, $\beta_j$ and $\gamma_j$, one
can construct $\mu_{1,2}$ explicitly, since $P_2$ is only quadratic in $\mu$.

\subsubsection{RLCs for the Enright methods}\label{section22}

\begin{figure}[H]
\centerline{\includegraphics[width=.68\textwidth]{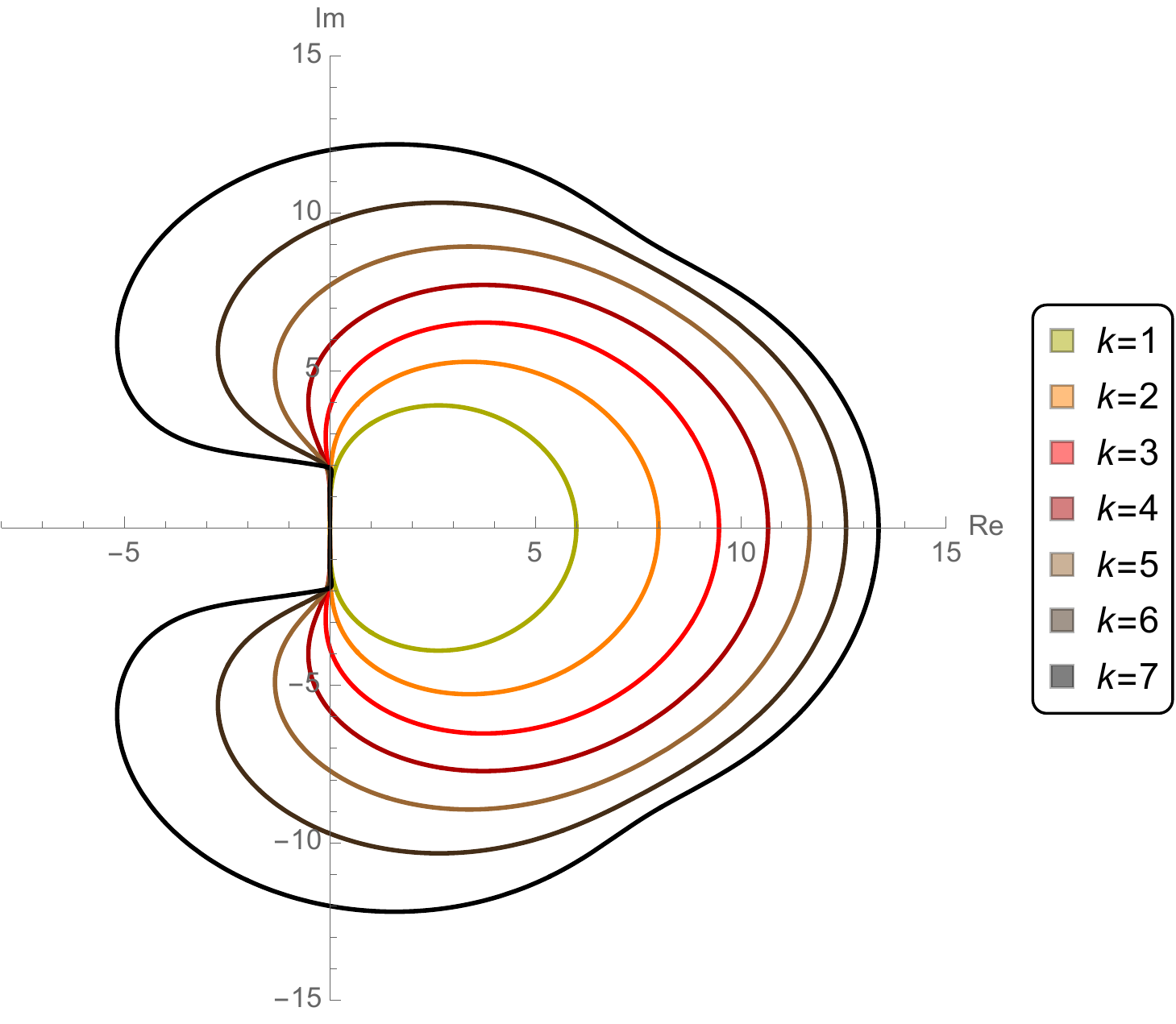}}
\caption{RLCs for the $k$-step Enright methods for $1\le k\le 7$. The stability
region of the method in each case is the unbounded component of $\Cc$.}
\label{fig:3}      
\end{figure}

\begin{figure}
\centerline{\includegraphics[width=.57\textwidth]{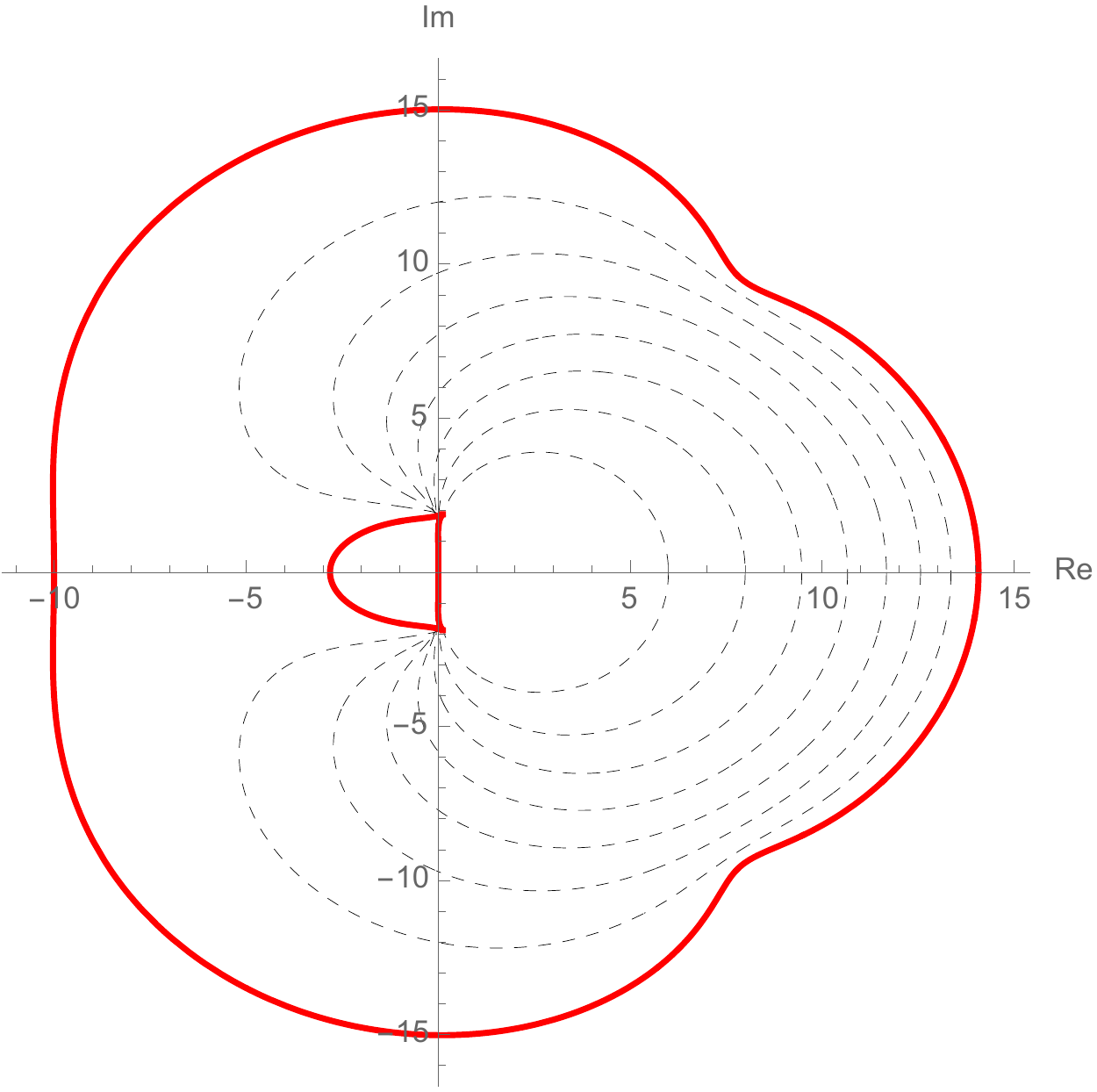}}
\caption{RLCs for the unstable $8$-step Enright method in red. The stability
region $\Ss$ is not connected, $\Cc\setminus\Ss$ is the annulus-like region. For comparison, 
the curves from Figure \ref{fig:3} are displayed as dashed gray curves.}
\label{fig:4}      
\end{figure}

The Enright methods are special cases of \eqref{2nddermultistepdefn}, and for $k\ge 1$ they are
defined \cite{hairerwanner} as
\begin{equation}\label{Enrightform}
y_{n+1}=y_n+h f_{n+1}-h\sum_{j=1}^k\left(\frac{1}{j}\left(\sum_{\ell=j}^k\nu_\ell\right)\nabla^j f_{n+1} \right)+h^2 \left(\sum_{\ell=0}^k \nu_\ell\right) g_{n+1}, 
\end{equation}
where \[\nu_\ell:=(-1)^\ell \int_0^1 (\tau-1)\binom{1-\tau}{\ell} d\tau 
\quad\quad  (0\le \ell \le k)\]
with the usual extension of the binomial coefficients.
From \eqref{Enrightform} one obtains the RLCs of
the Enright methods, see Figures \ref{fig:3} and \ref{fig:4}. The order of the $k$-step 
Enright method is $k+2$.

\section{Optimal sector inclusions}\label{determination}



\subsection{The RLC in implicit algebraic form}\label{section23}

Computing the stability
angle of a method with stability region $\Ss$ is equivalent to finding the slope of the unique line $L$ that passes 
through the origin, touches $\partial \Ss$ at some point in the open upper left half-plane such that 
$\partial \Ss$ lies on the right-hand side of $L$ (viewed from the origin) in this quadrant. This last requirement is
necessary since $\partial \Ss\cap L$ can consist of more points, even in the open upper left half-plane,
see Figure \ref{splitfigure3}. 

Assume now that $\partial \Ss$ can be represented by the RLC of the method (cf.~Remark \ref{nonvanishingremark}).  As we have seen, 
the RLC
is the image of the function $\mu(\cdot)$ in \eqref{muasaratio} for LMMs, or the union of the images of the functions $\mu_{1,2}(\cdot)$ in \eqref{mu12asaquadratic} for 
second-derivative multistep methods. The function $\mu$ is given
as a simple ratio, but to get the explicit forms of $\mu_{1,2}$, 
one should solve a quadratic equation. As the value of $k$ gets larger, these explicit formulae for 
$\mu_{1,2}$ 
 corresponding to a $k$-step second-derivative multistep method
become more and more complicated. Moreover, obtaining explicit and practically useful 
parametrized formulae for the
RLCs associated with multistep methods based on higher-than-second order derivatives
would almost be impossible.

To avoid these difficulties, we now describe a more general and effective technique which reduces the determination of the stability angles to the solution of a suitable system of polynomial equations. Let us consider 
the equation $\Phi(e^{i\vartheta},\mu)=0$ (see \eqref{phidef}).
By using the
well-known Weierstrass substitution \cite[pp.~382-383]{spivak}
\[
\vartheta = 2\arctan(t)\quad\quad (t\in\mathbb{R}),
\]
we have
$e^{i\vartheta}=
{(i-t)}/{(i+t)}$;
so instead of solving $\Phi(e^{i\vartheta},\mu)=0$ for $\mu$, we can solve 
\begin{equation}\label{Pinrationalfunctionform}
\Phi\left(\frac{i-t}{i+t}, \mu\right)=0 
\end{equation}
without trigonometric functions. Notice 
that originally we have $\vartheta\in [0,2\pi]$ in $e^{i\vartheta}$, or equivalently, $\vartheta\in (-\pi,\pi]$, but $\pi$ is not in the range of the function $2\arctan$, therefore we define 
\[
M_{-1}:=\{ \mu\in \mathbb{C} : \Phi\left(e^{i\pi},\mu\right)=0\}
\]
to restore the missing $\mu$ value(s) due to the reparametrization. 
Then, clearly, \eqref{Pinrationalfunctionform} can be brought to the form
$Q(t,\mu)/R(t)=0$ with some (complex) polynomials $Q$ and $R$. By writing
$\mu=a+b i$ ($a, b \in\mathbb{R}$) we get that 
there exist two real polynomials $Q_\text{re}:\mathbb{R}^3\to\mathbb{R}$ and 
$Q_\text{im}:\mathbb{R}^3\to\mathbb{R}$ such that the solutions
of $Q(t,\mu)=0$ for any fixed $t\in\mathbb{R}$ are obtained as the solutions of the system
\begin{equation}
\left\{
\begin{aligned}\label{Qreimsystem}
    Q_\text{re}(t,a,b)  &= 0 \quad\\ 
    Q_\text{im}(t,a,b)  &= 0.
\end{aligned}
\right.
\end{equation}
Now we eliminate $t$ by taking the resultant \cite{gelfand}
of $Q_\text{re}$ and $Q_\text{im}$ with respect to this parameter, and get that there exists a real
polynomial $F:\mathbb{R}^2\to\mathbb{R}$ such that if \eqref{Qreimsystem} holds
for some $t\in\mathbb{R}$, then $F(a,b)=0$ should hold with some $a, b \in \mathbb{R}$.
Hence, after identifying $\mathbb{C}$ with $\mathbb{R}^2$, we see that
the RLC can be represented as the implicit algebraic curve $C\cup M_{-1}$ with
$C:=\{ (a,b)\in \mathbb{R}^2 : F(a,b)=0\}$. Assuming that the set $M_{-1}$ is finite (it has at most two elements in the case of the BDF and Enright methods we are interested in), we ignore this component and focus only on $C$.
Suppose now that a line $L$ passes
through the origin and touches 
$C$ in the open 
upper left half-plane at some $(a_0, b_0)$ with $a_0<0<b_0$.
By assuming that
$C$ can be represented locally as the graph of an implicit function
near $(a_0,b_0)\in C$, we easily get,
by differentiating $a\mapsto F(a,b(a))$, 
that $(a_0, b_0)$ satisfies
\begin{equation}\label{FpartialFsystem}
\left\{
\begin{aligned}
    F(a_0,b_0) &= 0 \quad\\ 
    a_0\cdot\partial_1 F(a_0,b_0)+b_0 \cdot\partial_2 F(a_0,b_0)  &= 0\\
    a_0 & < 0\\
    b_0  & > 0.
\end{aligned}
\right.
\end{equation}
By taking again the resultant of the first two polynomial equations, one of the variables, say $b_0$, is 
eliminated. The resulting univariate polynomial yields in the general case finitely many possible $a_0$ values to choose 
from. With $\al$ denoting the angle (in radians) between $L$ and the negative half of the real axis, we get
that $\tan(\al)=-{b_0}/{a_0}$. To select the appropriate solution
$(a_0,b_0)$  (and hence the appropriate tangent line $L$), we verify in the concrete case that $(a_0,b_0)\in \partial \Ss\subset\mathbb{C}=\mathbb{R}^2$, and determine whether $\partial \Ss$ lies on the right-hand side of $L$. The appropriately chosen $\al$ angle then yields the
desired stability angle.


\subsection{Results for the BDF methods}\label{BDFsection}

\begin{figure}[H]
\subfigure{
\includegraphics[width=0.47\textwidth]{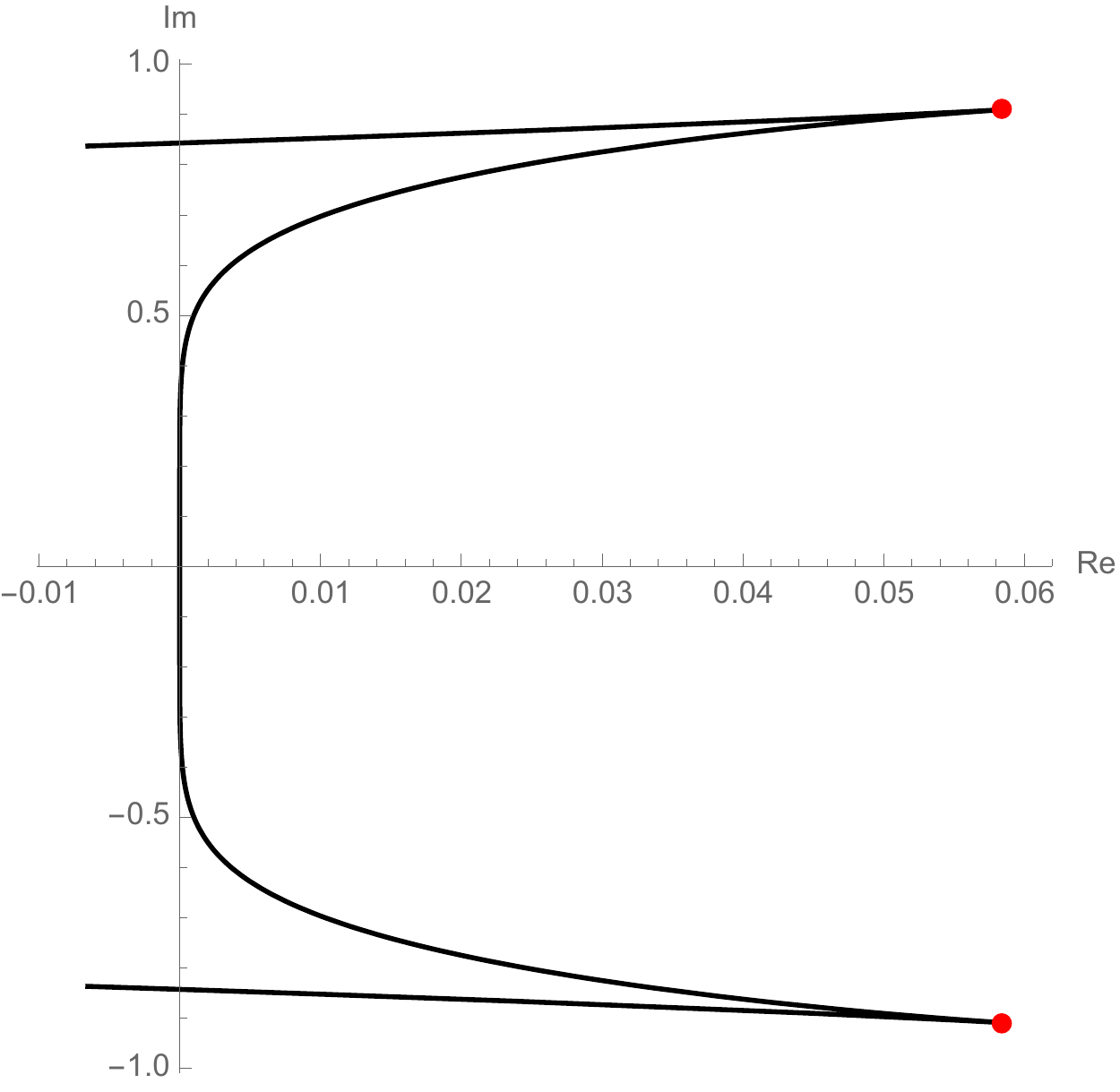}}
\subfigure{
\includegraphics[width=0.47\textwidth]{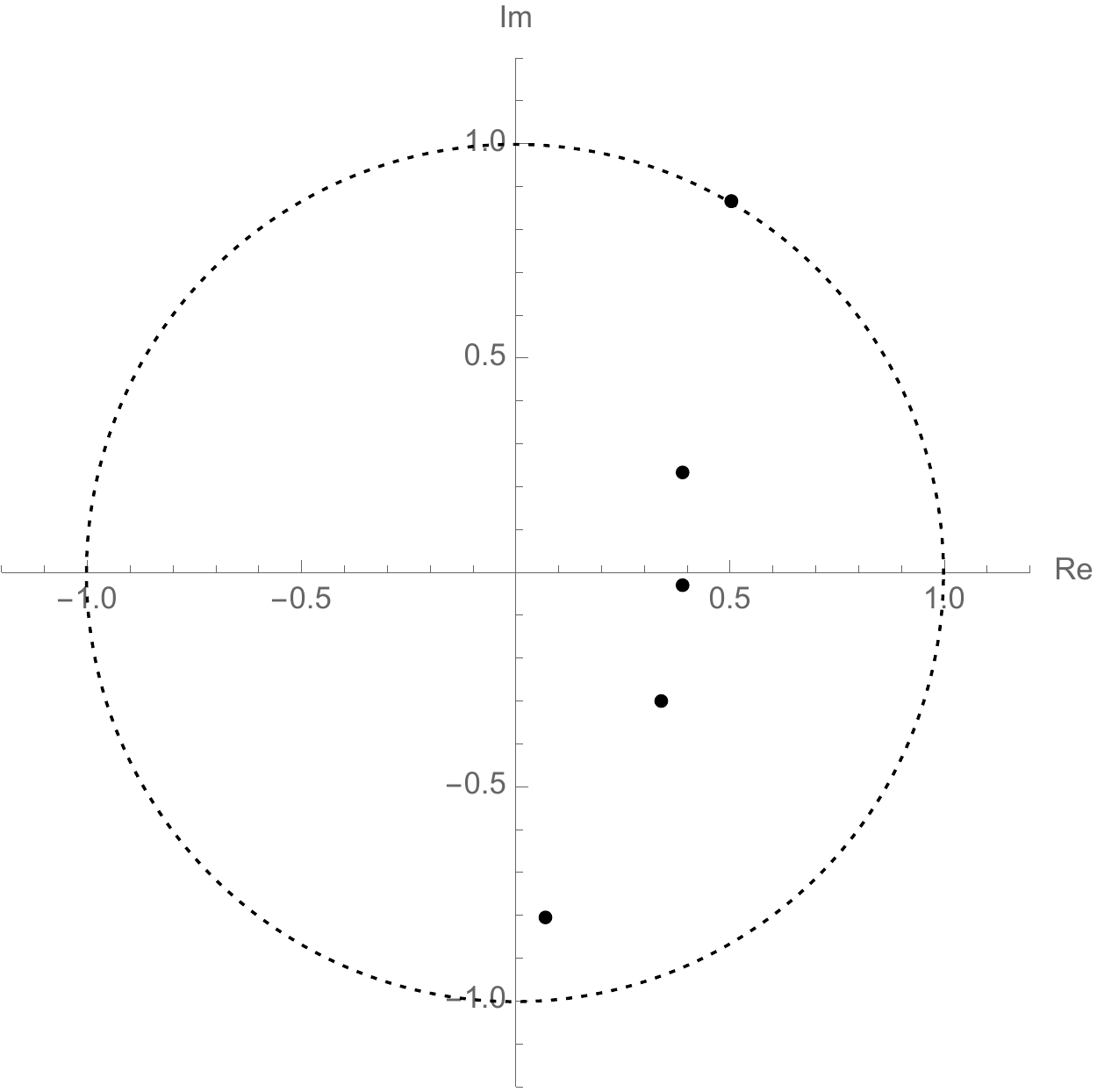}}
\caption{Cusp singularities of $\partial \Ss$ for the $6$-step BDF method, denoted
by red dots in the left figure. The singularities
are located at $\mu_{\pm}:=\frac{7}{120}\pm i\, \frac{21 \sqrt{3}}{40}\approx 0.0583\pm  
0.9093 i$. For each such $\mu$ value,
 $P_1(\cdot,\mu)$ in \eqref{P1def} has a double root with modulus
equal to 1. Therefore $\mu_\pm\in\partial\Ss\setminus \Ss$, 
hence this $\Ss$ is not closed. The right figure depicts the 6 roots of 
$P_1(\cdot,\mu_+)$, and the double root is located at $\frac{1}{2} \left(1+i \sqrt{3}\right)$ (note that $\mu_+\in\mathbb{C}\setminus\mathbb{R}$, so
these roots are not symmetric with respect to the real axis).
\label{BDF_6_cusp}}
\end{figure}

The simplest non-trivial case illustrating the steps in Section \ref{section23} is 
the determination of the stability angle for the $3$-step BDF method.
Formula \eqref{BDFconcreteformula} with $k=3$ yields the following trigonometric parametrization of the RLC in $\mathbb{R}^2$ after a simplification:
\[
[0,2\pi]\ni\vartheta\mapsto\mu(\vartheta):=
\]
\[\left(\frac{4}{3} \sin ^4\left(\frac{\vartheta }{2}\right) (1-4 \cos (\vartheta )),
\frac{\sin (\vartheta )}{3}  \left[2 (\cos (2 \vartheta )+5)-9 \cos (\vartheta )\right]\right).
\]
After eliminating the trigonometric functions, \eqref{Pinrationalfunctionform} can be written as
\[
\frac{Q(t,\mu)}{R(t)}=\frac{3 \mu  t^3-20 t^3-9 \mu  t+6 t+i \left(3 \mu -9 \mu  t^2+18 t^2\right)}{3 (t-i)^3}=0.
\]
Then $Q_\text{re}$ and $Q_\text{im}$ in \eqref{Qreimsystem} become
\[
\left\{
\begin{aligned}
    3 a t^3-9 a t+9 b t^2-3 b-20 t^3+6 t &= 0 \quad\\ 
    -9 a t^2+3 a+3 b t^3-9 b t+18 t^2  &= 0.
\end{aligned}
\right.
\]
We eliminate $t$ from this system and obtain
\[
F(a,b):=432 \left[108 a^6-1188 a^5+9 a^4 \left(36 b^2+439\right)-2 a^3 \left(1188 b^2+3121\right)\right.+\]\[\left.
9 a^2 \left(36 b^4+394 b^2+547\right)-54 a \left(22 b^4+17 b^2+30\right)+27
   b^4 \left(4 b^2-15\right)\right].
\]
Now $b$ is eliminated from the first two equations of \eqref{FpartialFsystem}, and we get that the
possible choices for $a_0$ are the negative real roots of
\[
a^4 (24 a-25)^4 (5324 a+405)^2 \left(6 a^2-13 a+9\right)^2=0,
\]
yielding the unique value $a_0=-{405}/{5324}$. Substituting this $a_0$ into \eqref{FpartialFsystem} 
we get the unique value $b_0={987 \sqrt{35}}/{5324}$, hence
$\tan(\al)=-b_0/a_0=({329 \sqrt{7/5}})/{27}$ is the only possible value for the tangent of the
stability angle. Finally, we verify that the corresponding tangent line $L$ passing through
the origin has no other intersection
point with $\partial \Ss$ in the open upper left quadrant, and $\partial \Ss$ lies on the right
side of $L$.
\begin{remark} The above RLC for the $3$-step BDF method can also be parametrized 
 as
\[
\mathbb{R}\ni t\mapsto \left(\frac{4 t^4 \left(5 t^2-3\right)}{3 \left(t^2+1\right)^3},\frac{2 t \left(21 t^4+8 t^2+3\right)}{3 \left(t^2+1\right)^3}\right)\in\mathbb{R}^2.
\]
Here $M_{-1}=\{(20/3, 0)\}\subset\mathbb{R}^2$, corresponding to the $t\to \pm \infty$ limiting value of the
parametrization.
\end{remark}
The remaining stability angle values for $4\le k\le 6$ can be computed analogously, so Table \ref{tab:1}
shows only the final exact results. 
\begin{table}[h]
\caption{The exact stability angles $\alpha_k^\text{BDF}=\frac{180}{\pi}\arctan\left(c_k^\text{BDF}\right)^\circ$ of the BDF methods}
\label{tab:1}      
\centerline{\begin{tabular}{l|l|l}
\hline\noalign{\smallskip}
$k$ & $c_k^\text{BDF}$ & Approximate value of $\al_k^\text{BDF}$\\
\noalign{\smallskip}\hline\noalign{\smallskip}
3 & $\frac{329 \sqrt{\frac{7}{5}}}{27}$ & $86.032366860211647332^\circ$ \\
4 & $\frac{699 \sqrt{\frac{3}{2}}}{256}$ & $73.351670474578482110^\circ$ \\
5 & $\frac{1326107429}{25} \sqrt{\frac{62}{53860574450525125+1194498034900685 \sqrt{2033}}}$ & $51.839755836049910391^\circ$ \\
6 & $\frac{45503}{10125 \sqrt{195}}$ & $17.839777792245700101^\circ$ \\
\noalign{\smallskip}\hline
\end{tabular}}
\end{table}
\begin{remark}
For $3\le k\le 6$, the BDF stability region includes an interval along the imaginary axis and 
containing the origin if and only if $k=5$ or $k=6$. For $k=5$ and $k=6$ the 
two intervals are
\[
\{z\in\mathbb{C} : \mathrm{Re}(z)=0,\, |\mathrm{Im}(z)|\le \frac{1}{12\sqrt{2}} \sqrt{12775-387 \sqrt{1065}}\approx 0.710\}\subset \Ss
\]
and 
\[
\{z\in\mathbb{C} : \mathrm{Re}(z)=0,\, |\mathrm{Im}(z)|\le \frac{7}{20} \sqrt{1263-336 \sqrt{14}}\approx 0.843\}\subset \Ss,
\]
respectively.
\end{remark}
\begin{remark}
The boundary curve of the stability region of the $6$-step BDF method
contains two cusp singularities, see Figure \ref{BDF_6_cusp} (and compare with 
Figure \ref{trueboundary}). No other $\partial \Ss$ curve
has this type of degeneracy in the BDF family for $1\le k\le 5$ or $k=7$. 
Since the cusp points for $k=6$ are not part of $\Ss$, the stability region
in this case is not closed (nor open).
\end{remark}

\subsection{Results for the Enright methods}\label{Enrightsection}

By applying
the algorithm described in Section \ref{section23}, we can exactly determine the 
stability angles for the Enright methods, see Table \ref{tab:2}.
But since the $c_k^E$ values are much more complicated algebraic numbers
than the 
corresponding $c_k^B$ constants in Table \ref{tab:1}, 
Table \ref{tab:2} contains only a numerical approximation to the exact
stability angles.
\begin{table}[H]
\caption{Stability angles $\alpha_k^\text{Enr}=\frac{180}{\pi}\arctan\left(c_k^\text{Enr}\right)^\circ$ of the Enright methods}
\label{tab:2}      
\centerline{\begin{tabular}{l|l|l}
\hline\noalign{\smallskip}
$k$ & Approximate value of $c_k^\text{Enr}$ & Approximate value of $\alpha_k^\text{Enr}$  \\
\noalign{\smallskip}\hline\noalign{\smallskip}
3 & 27.056933440109472532101963 & $87.8833627693413031369003498^\circ$ \\
4 & 7.1406622283653916403051061 & $82.0279713768712835947479188^\circ$ \\
5 & 3.2907685080317853840110455 & $73.0970020659749082763655203^\circ$ \\
6 & 1.7285146253131256601603521 & $59.9492702555400766770433070^\circ$ \\
7 & 0.7703217281441388675578954 & $37.6078417405752150238159031^\circ$ \\
\noalign{\smallskip}\hline
\end{tabular}}
\end{table}

\begin{remark}
By rounding the values of $\al_k^\text{Enr}$  
given in Table \ref{tab:2} to two decimal places, we recover the approximate values of these stability
angles in \cite[Chapter V.3, Table 3.1]{hairerwanner}.
\end{remark}

It turns out that $c_3^\text{Enr}$ is an algebraic number of degree 22,
being the unique positive root of the following even polynomial with coefficients
\[
\{6621625501626720011970719022734459520000000000000000, 0,\]\[4744945665370497147850526235135397935643117766707200000, 0,\]\[74537179754361052063480563770102869789636567887828480000, 0,\]\[417809113212221868517393954677075422852686053100794277975, 0,\]\[1103592881533264097533512931940128409045933472020943607320, 0,\]\[1780216754145335084531442707748395556646595339402356863603, 0,\]\[2028417751642933570985301304414377204911584843581604760752, 0,\]\[1720629215811045658880293770988465046952673868659037700813, 0,\]\[1065257770963658030926145190690110109450795207237154063632, 0,\]\[451976742777053443392779380035051991794204051855298481913, 0,\]\[117280744006618927204325767614876515512652225395198902600, 0,\]\[14037302894263476230042573549418427869442188056651130000\}.
\]
\begin{remark}
Besides the stability angle, there are  
other measures of stability for $A(\al)$-stable methods. One of these characteristics
is the \emph{stiff stability abscissa}, being the smallest constant $D>0$  such that
$\{ z\in\mathbb{C} : \mathrm{Re}(z)\le -D\}\subset \Ss$. For example, for the $3$-step Enright method, 
Table 3.1 in \cite[Chapter V.3]{hairerwanner} contains the approximate value $D\approx 0.103$. By
using our implicit representation of $\partial \Ss$, it is straightforward to determine
the exact value of $D\approx 0.10341810907195$; it is an algebraic number of degree 12, and the total number of digits in the coefficients of its defining integer polynomial is 529.
\end{remark}

As for the $k=4$ case, the algebraic degree of $c_4^\text{Enr}$ is 28. The constants $c_5^\text{Enr}$, $c_6^\text{Enr}$ and $c_7^\text{Enr}$ can be given
as roots of increasingly more involved integer polynomials, so we do not
reproduce these
polynomials here. During the computations in the $k=7$ case, for example, 
we had to manipulate intermediate polynomials of degree of a few hundred, or polynomials
with a total number of coefficient digits of approximately 470000. We could describe the final defining polynomial for $c_7^\text{Enr}$ by $\approx 175000$
characters in \textit{Mathematica}.

\begin{figure}[H]
\subfigure{
\includegraphics[width=0.422\textwidth]{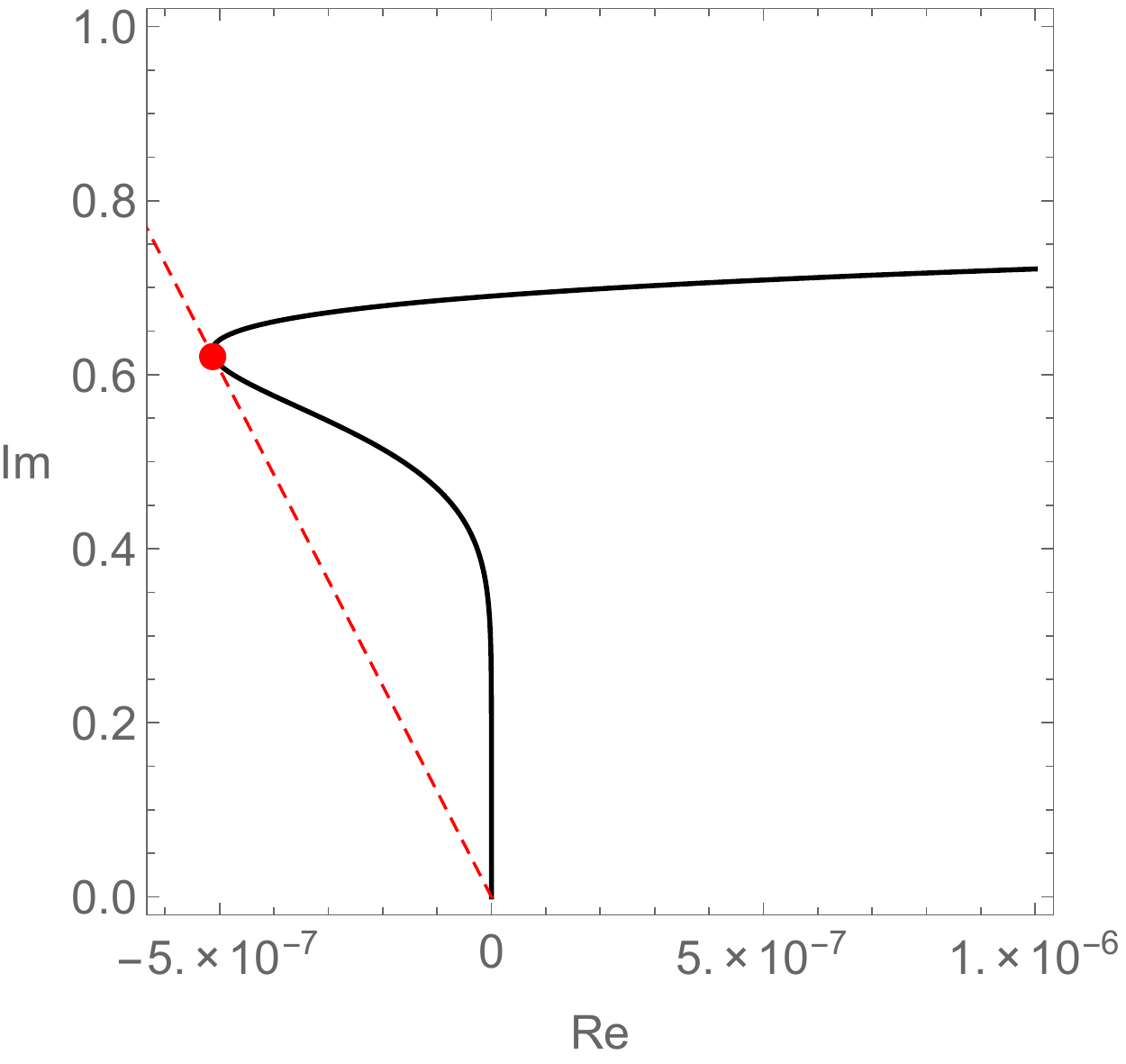}}
\subfigure{
\includegraphics[width=0.4\textwidth]{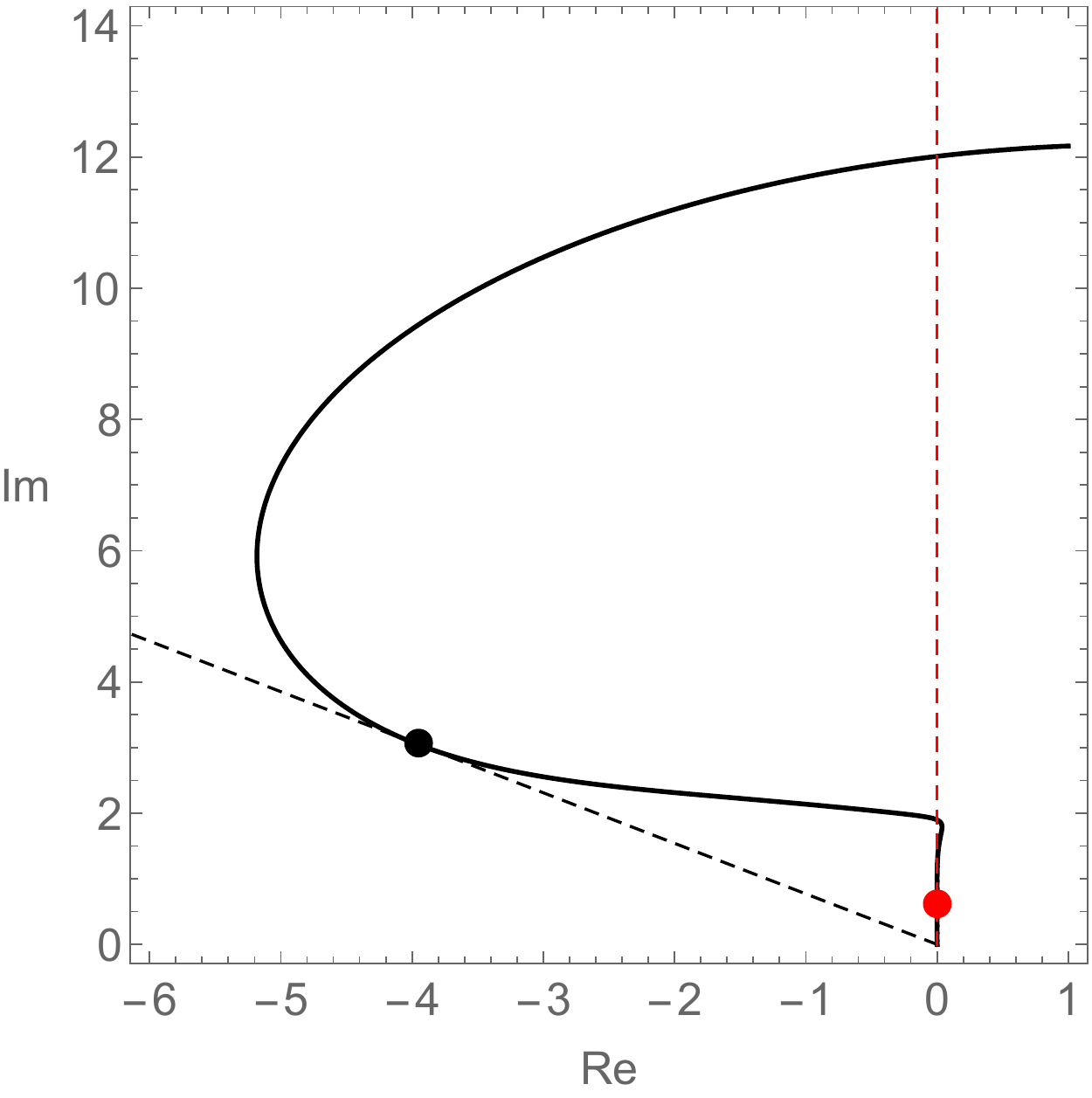}}
\caption{Part of the boundary of the stability region of the $7$-step Enright method near
the origin (solid black curve) together with the two (dashed red and black) lines that 
pass through the origin and are locally tangent to the boundary curve at some point in the open 
upper left half-plane. 
Due to the scaling, the
dashed black line is seen only in the larger plot window on the right. The stability angle 
$\al_7^E\approx 37.6^\circ$ of the method is determined by the
dashed black line; the red line has additional intersection points with the boundary curve.  The angle between the dashed red line and the negative half of the real axis has also been computed exactly; its approximate
value is $\approx 89.9999527^\circ$.
\label{splitfigure3}}
\end{figure}

\begin{remark}
Let us consider the Enright stability region corresponding to $k=7$. As we already remarked earlier,  
there are exactly two lines
that pass through the origin and are locally tangent to the boundary curve at some point in the open 
upper left half-plane, see Figure \ref{splitfigure3}. Within the BDF family for $1\le k\le 6$ or
in the Enright family for $1\le k\le 7$, this phenomenon occurs only in the present case.
\end{remark}

\section{Optimal disk inclusions}\label{diskinclusionsection}

\begin{figure}
\subfigure{
\includegraphics[width=0.45\textwidth]{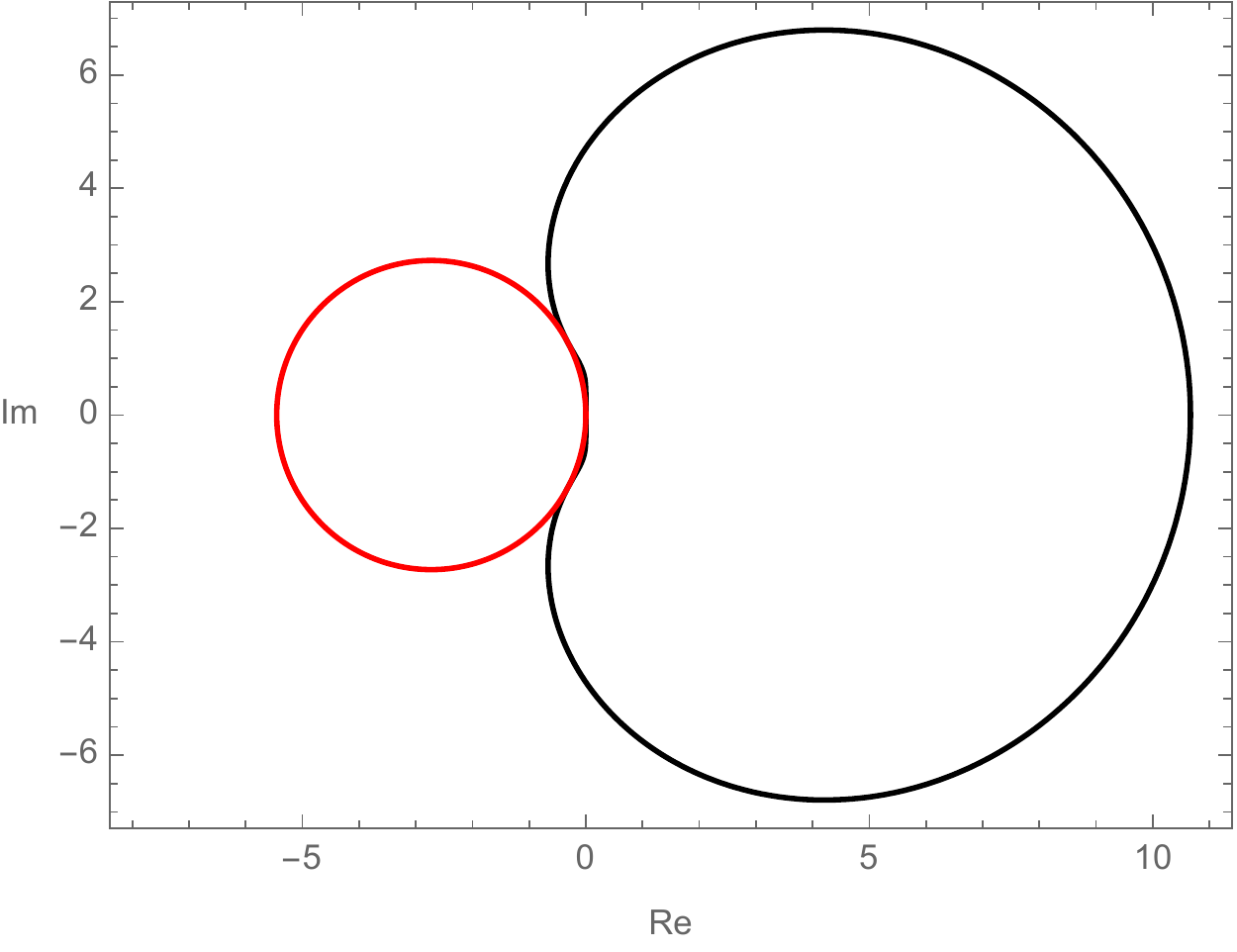}}
\subfigure{
\includegraphics[width=0.45\textwidth]{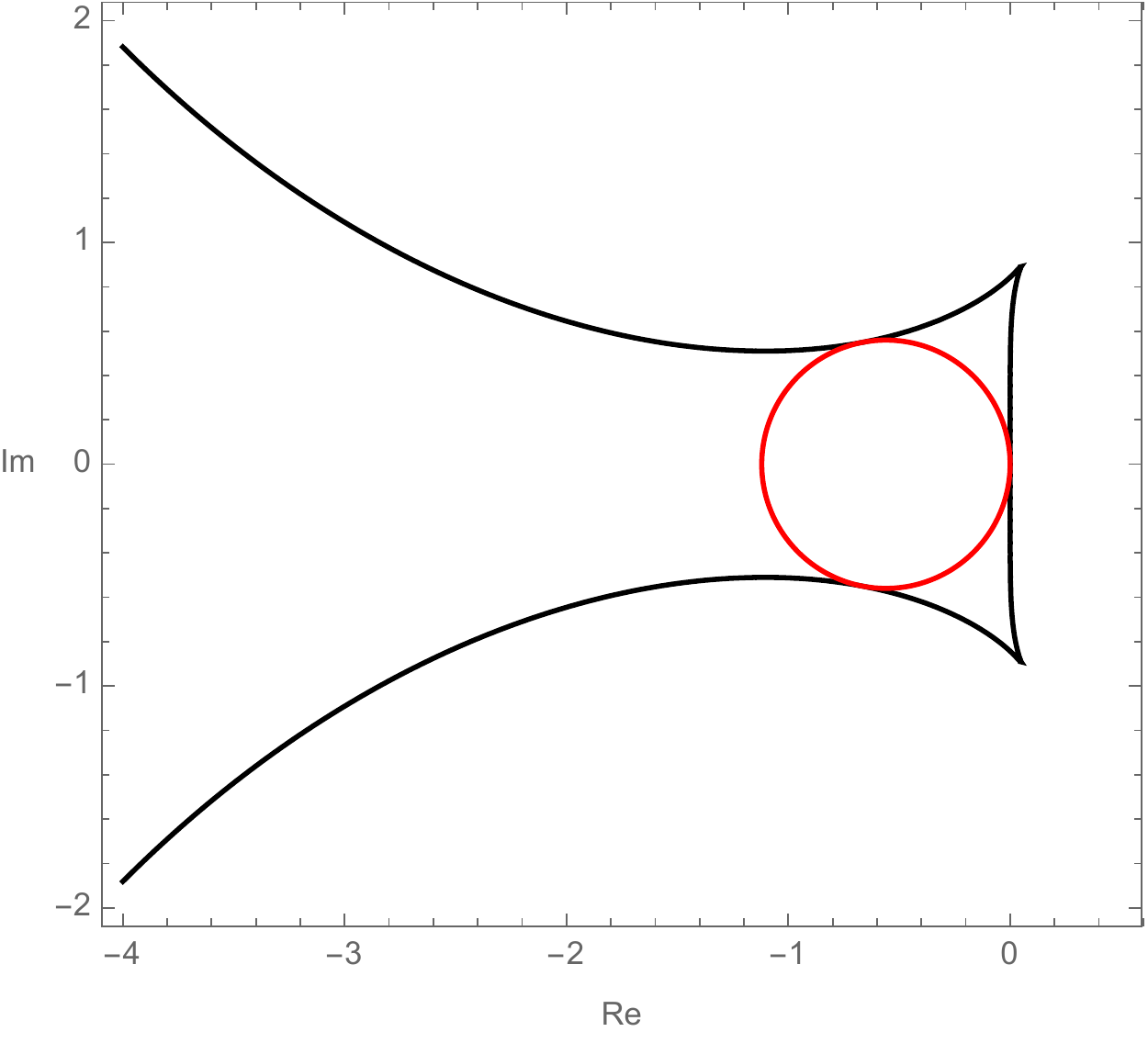}}
\caption{The largest inscribed disk $|z+r|\le r$ (with red boundary)  in the stability region of the $k$-step BDF method for $k=4$ (left) and $k=6$ (right), see Table \ref{tab:3}.}
\label{BDF6maxdisk}      
\end{figure}

As for the largest inscribed disk in the stability region $\Ss$, we again 
expect---similarly to Section \ref{section23}---that $\partial \Ss$ (or the RLC) and the optimal disk possess a common tangent line (with point of tangency different from the origin). 
By using
\begin{itemize}
\item  the implicit algebraic form $F(a,b)=0$ of the RLC,
\item  the implicit equation $(a+r)^2+b^2-r^2=0$ for the boundary of the inscribed disk,
\item  and the condition for a common tangent line
\[-\frac{\partial_a F(a,b)}{\partial_b F(a,b)}=-\frac{\partial_a \Big((a+r)^2+b^2-r^2\Big)}{\partial_b \Big((a+r)^2+b^2-r^2\Big)},\]
\end{itemize}
we obtain a system of 3 polynomial equations in 3 unknowns $(a,b,r)$. By taking resultants and successively eliminating the variables ($a, b$), we obtain a univariate polynomial in $r$ whose positive
root will yield the optimum value of the stability radius.
The exact optimal stability radii $r_k^\text{\,BDF}$ for the $k$-step BDF methods ($3\le k \le 6$) are
found in Table \ref{tab:3}; see Figure \ref{BDF6maxdisk} also.
The degree of the algebraic number $r_k^\text{\,BDF}$ is $2, 3, 5, 5$ for $3\le k \le 6$, respectively.
\begin{remark} It is quite surprising that the algebraic numbers listed in Table \ref{tab:3} have such a low degree for the following reasons.   
For the 3-step BDF method, the univariate polynomial $r$ mentioned above
has degree 28, but it can be split into several factors of lower degree, and has a unique positive root $r_3^\text{\,BDF}\approx 7.0497$. For the 4-step BDF method, the corresponding $r$-polynomial has degree 52 and a unique positive root $\approx 2.7272$.  The 
 $r$-polynomial for the 5-step BDF method has degree 88 and a unique positive root $\approx 1.3579$. 
Finally, the $r$-polynomial for the 6-step BDF method has degree 128, and a unique positive
root $\approx 0.5599$. 
\end{remark}

\begin{table}[h]
\caption{The exact stability radii $r_k^\text{\,BDF}$ of the BDF methods}
\label{tab:3}       
\centerline{\begin{tabular}{l|l|l}
\hline\noalign{\smallskip}
$k$ & $r_3^\text{\,BDF}$ is equal to / $r_{4,5,6}^\text{\,BDF}$ is a root of the polynomial & Approximate value of $r_k^\text{\,BDF}$\\
\noalign{\smallskip}\hline\noalign{\smallskip}
3 & $\left(17+8 \sqrt{10}\right)/6$ & $7.049703546891172$ \\
4 & $\{18432, 2172, -100855, -114975\}$ & $2.727199466336645$ \\
5 & $\{2944512000, 260854387200, 679386763440,$ &  \\
 & $ \ 266052478296, -1280160594125, -1354065829875\}$ & $1.357947301777465$ \\
6 & $\{141717600000, 558150393600, 1112790780640,$ &  \\
 & $\ 948530730784, -119637602525, -488414721375\}$ & $0.559931687924882$ \\
\noalign{\smallskip}\hline
\end{tabular}}
\end{table}

\section{Optimal stability angle in a family of multistep methods}\label{section7}

In \cite{hundsdorferruuth}, ODEs of the form $u'(t)=F(u(t))+G(u(t))$, $u(0)=u_0$
are considered, with $F$ and $G$ representing 
non-stiff and stiff parts of the equation, respectively.
To solve these equations numerically, the authors construct 
several implicit-explicit (IMEX) LMMs, and thoroughly analyze them 
from the viewpoint of numerical monotonicity, boundedness and stability. Their analysis
involves finding optimal methods with respect to various criteria in certain families. 

Here we take their simplest case study from \cite[Section 3.2.1]{hundsdorferruuth},
a 2nd-order, 3-step explicit method augmented by an implicit method
(note that we changed their notation from $b_j$ to $\beta_j$):
\begin{equation}\label{modelexample}
u_n=\frac{3}{4}u_{n-1}+\frac{1}{4}u_{n-3}+\frac{3}{2}\Delta t \cdot F_{n-1}+
\sum_{j=0}^3\beta_j \Delta t\cdot G_{n-j}.
\end{equation}
The values of $\beta_2:=-3 \beta _0-2 \beta _1+3$ and $\beta_3:=2 \beta _0+\beta _1-\frac{3}{2}$
are determined from the order conditions, so \eqref{modelexample} becomes a 2-parameter family of methods, with real parameters $\beta_1$ and $\beta_0$. 
The three figures in \cite[Figure 1]{hundsdorferruuth} then depict the 
$A(\al)$-stability angles, the ``damping factors'' and the ``absolute error constants'',
respectively, of members of the family \eqref{modelexample}. In what follows, we 
do not consider these last two categories but focus only on the leftmost figure in
 \cite[Figure 1]{hundsdorferruuth}---as the authors conclude in \cite[Section 3.2.1]{hundsdorferruuth}, 
a method with large stability angle does not necessarily have a good damping factor or
a small error constant, and vice versa; the different optimization criteria are often conflicting. In other
words, our goal in this section is to find the IMEX method in the family \eqref{modelexample} with the largest stability angle. 

To begin the $A(\al)$-stability investigation, the authors of \cite{hundsdorferruuth} 
define the usual linear test functions $F(u):=\hat{\lambda} u$ and
$G(u):=\lambda u$. They then assume that $\Delta t\cdot \hat{\lambda}=i\eta$ and 
$\Delta t\cdot \lambda =\xi$ with $\eta\in\mathbb{R}$ and 
$\mathbb{R}\ni \xi\le 0$: this choice is relevant
``for example, for advection-diffusion equations if central finite differences or spectral
approximations are used in space''. These assumptions lead to the following 
characteristic polynomial of the IMEX multistep family, see \cite[(2.4)--(2.7)]{hundsdorferruuth}:
\begin{equation}\label{originalcharpoly}
\mathbb{C}\ni\zeta\mapsto\zeta^3-\left(\frac{3}{4} \zeta ^{2}+\frac{1}{4}\right)-i\eta \left(\frac{3}{2} \zeta ^2\right) -\xi  \left(\sum _{j=0}^3 \beta _j \zeta ^{3-j}\right).
\end{equation}
To create the leftmost figure in \cite[Figure 1]{hundsdorferruuth} \textit{approximately} indicating the optimal stability angle within the family, the authors use \eqref{originalcharpoly} to construct the RLCs and study these curves ``for $\xi\to-\infty$'' to \textit{estimate} the stability 
angles\footnote{\label{footnote2}When 
the 
stability angle $\al$ of a method is defined in \cite[Sections 2.3 and 3.2.1]{hundsdorferruuth} notice that 
we should require that the sector 
\[
\xi\le 0, \quad |\eta/\xi|\le \tan(\al)\quad \text{ with angle } \al\le\pi/2
\]
be included in the stability region in the $(\xi,\eta)$-plane
(with the $\xi=0$ and $\al=\pi/2$ cases interpreted appropriately). In other
words, $\arctan(\al)$ in \cite{hundsdorferruuth} is to be replaced by $\tan(\al)$, otherwise
the sector would not ``open wide enough''
and $A$-stability would not be recovered in the $\al\to \pi/2^-$ limit.
See also Footnote \ref{footnote3}.}.

In the rest of this section we confirm their numerical findings,
 but we
solve the optimization problem \textit{rigorously and exactly}. We have selected this family 
\eqref{modelexample} because the final 
result---the optimal stability angle---has a particularly simple 
form (see our Theorem \ref{wedgethm} below), and, at the same time, our straightforward approach based on the theorems cited in Section \ref{schurcohnsection}  
is readily illustrated. We emphasize that our analysis avoids the construction of the 
RLCs: as we have seen (for example, in Figure \ref{trueboundary}), they may have complicated
self-intersections, and it is often not obvious a priori 
whether a particular segment of the RLC coincides with the stability region boundary or not.

\subsection{Summary of the main steps and results}

By rearranging \eqref{originalcharpoly} and inserting the values of $\beta_2$ and $\beta_3$
given below \eqref{modelexample}, we define 
\[
P_{\beta_1,\beta_0}(\zeta,\xi,\eta):=\left(1-\beta _0 \xi \right)\zeta ^3 -
 \left(\frac{3}{4}+\beta _1 \xi +\frac{3 i \eta }{2}\right)\zeta ^2+\]
\begin{equation}\label{Pbetadef}
\xi\left(3 \beta _0+2 \beta _1-3\right) \zeta-
   \left(\frac{1}{4}+2 \beta _0 \xi + \beta _1 \xi -\frac{3}{2} \xi \right),
\end{equation}
where $\zeta\in\mathbb{C}$, $(\beta_1, \beta_0)\in\mathbb{R}^2$, $\xi\le 0$ and
$\eta\in\mathbb{R}$. Our goal is to find the parameters $(\beta_1,\beta_0)$   such that
the stability region 
\begin{equation}\label{Sbetadef}
\Ss_{\beta_1,\beta_0}:=\{ (\xi, \eta) \in\mathbb{R}^2: \xi\le 0,\, \eta\in\mathbb{R}, \,
P_{\beta_1,\beta_0}(\cdot,\xi,\eta)\in\svn\}
\end{equation}
contains the infinite sector 
\[
{\mathcal{A}}_m:=\{(\xi,\eta) \in\mathbb{R}^2 : \xi\le 0, \, \eta\in\mathbb{R}, \, |\eta|\le m |\xi|\}
\]
with the largest $m>0$ in the definition of $A(\al)$-stability. In other words, we are to find $(\beta_1,\beta_0)$ 
such that 
\begin{equation}\label{inclusion}
{\mathcal{A}}_m\subset \Ss_{\beta_1,\beta_0}
\end{equation}
holds with the largest possible $m>0$. Note that for convenience we have identified 
$\mathbb{C}$ with $\mathbb{R}^2$, hence stability regions in this section are subsets
of $\mathbb{R}^2$.

As a first step, Lemma \ref{Wlemma} below yields a necessary condition for the inclusion \eqref{inclusion}. In its proof---presented in 
Appendix \ref{sectionsomeneccond}---we use the argument proposed in \cite{hundsdorferruuth} and consider
the $\xi\to-\infty$, $\eta=0$
limiting values.
At this point it is convenient to recall the notion
of $\overset{\circ}{\text{A}}$-stability \cite[Chapter V.2]{hairerwanner}: 
a method is $\overset{\circ}{\text{A}}$-stable, if its stability region includes the non-positive reals
$\{\xi\in \mathbb{R} : \xi\le 0\}$. 
 Clearly, 
\[
A(\alpha)\text{-stability with some } \alpha>0 \implies  
\overset{\circ}{\text{A}}\text{-stability.}
\]

\begin{figure}[H]
\centerline{\includegraphics[width=.58\textwidth]{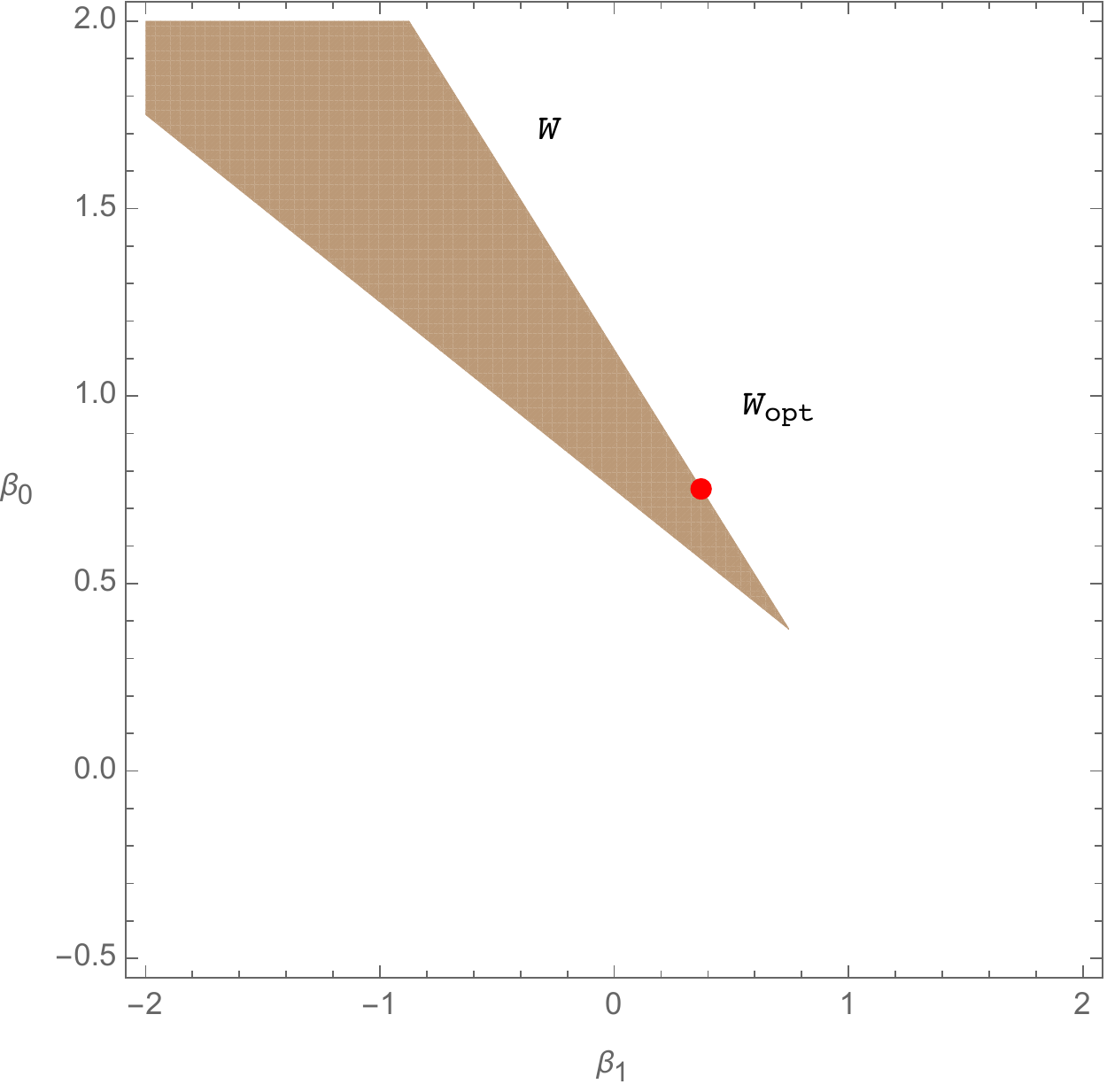}}
\caption{The figure shows the set $W$ defined in \eqref{Wjdef} and the optimal
parameter choice $W_{\text{opt}}$ with $(\beta_1,\beta_0)=(3/8, 3/4)$ determined in 
Appendix \ref{optimalmethodW2}. Interestingly,
the coordinates of the vertex of the wedge $W$ are $(3/4, 3/8)$.}
\label{fig:wedge}      
\end{figure}

\begin{lemma}\label{Wlemma} Let us define 
\begin{equation}\label{Wjdef}
W:=\left\{(\beta_1,\beta_0) : \beta _1\leq \frac{3}{4}, \ \frac{3-2 \beta _1}{4}\leq \beta _0\leq \frac{9-8 \beta _1}{8}\right\}.
\end{equation}
Then a method of the form \eqref{modelexample}
is \emph{not $\overset{\circ}{\text{A}}$-stable} for $(\beta_1,\beta_0)\notin W$. 
\end{lemma}
As a consequence, 
from now on we can assume $(\beta_1,\beta_0)\in W$, 
see Figure \ref{fig:wedge}. Note that the orientation of the axes in Figure \ref{fig:wedge} and in the leftmost figure in \cite[Figure 1]{hundsdorferruuth} is the same: the $\beta_1$-axis is horizontal, while the $\beta_0$-axis is 
vertical. 
Lemma \ref{Wlemma} thus
also proves that the wedge-like object in the parameter space in the 
leftmost figure in \cite[Figure 1]{hundsdorferruuth} is indeed a perfect (infinite) wedge given by
$W$. 

\begin{remark}
The assumption $(\beta_1,\beta_0)\in W$ implies
$\beta_0>0$, so due to $\xi\le 0$, the leading coefficient of \eqref{Pbetadef}, $1-\beta _0 \xi$, 
cannot vanish (cf.~Remark \ref{nonvanishingremark}).
\end{remark}

Then in Appendix \ref{optimalmethodW2} we prove the main result of Section 
\ref{section7}.
\begin{theorem}\label{wedgethm}
Suppose that $(\beta_1,\beta_0)\in W$. Then the largest $m>0$ such that
\eqref{inclusion} holds is $m\equiv m_{\text{opt}}:=1/2$.
\end{theorem}

In the proof we show that finding the optimal
$(\beta_1,\beta_0)\in W$ is equivalent to finding the largest positive real root of a suitable
polynomial in $m$ with coefficients depending on $\beta_1$ and $\beta_0$. 
We verify that this optimal root is located at $m_{\text{opt}}$, 
corresponding to the unique method with
$(\beta_1,\beta_0)=W_{\text{opt}}:=(3/8, 3/4)\in W$ and represented as a red dot 
in the parameter space in Figure \ref{fig:wedge}. 
The black curve in the left half-plane in Figure \ref{fig:optimalIMEXmethod} is the boundary 
of the optimal stability region, 
and the dashed red lines bound the largest inscribed infinite sector
${\mathcal{A}}_{1/2}$:
the optimal stability angle satisfies $\tan(\al)=m_{\text{opt}}$. 
As a conclusion, the highest value in the scale adjacent to the leftmost 
figure in \cite[Figure 1]{hundsdorferruuth} should be exactly $\al=\arctan(1/2)\approx 
0.463648$, that is, $\alpha\approx 26.5651^\circ$.


\begin{remark}\label{remarktouch}
Unlike in Section \ref{section8} (see Remark \ref{closingremark16}), the boundary of the optimal sector
${\mathcal{A}}_{1/2}$ does not touch (or intersect) the boundary of the optimal stability region
$\Ss_{3/8, 3/4}$ in the open left half-plane.
\end{remark}

\begin{remark}
In \cite[Section 3.2.1, (3.4)--(3.5)]{hundsdorferruuth}, the stability angles for 
two particular schemes from the family \eqref{modelexample} are also approximated. 
For the IMEX-Shu(3,2) scheme
\[
u_n=\frac{3}{4}u_{n-1}+\frac{1}{4}u_{n-3}+\frac{3}{2}\Delta t \cdot F_{n-1}+
\]
\[
\frac{4}{9}\Delta t\cdot G_{n}+
\frac{2}{3}\Delta t\cdot G_{n-1}+
\frac{1}{3}\Delta t\cdot G_{n-2}+
\frac{1}{18}\Delta t\cdot G_{n-3}
\]
they obtain $\alpha_{\emph{Shu}}\approx 0.06$, and for the IMEX-SG(3,2) scheme
\[
u_n=\frac{3}{4}u_{n-1}+\frac{1}{4}u_{n-3}+\frac{3}{2}\Delta t \cdot F_{n-1}+
\Delta t\cdot G_{n}+
\frac{1}{2}\Delta t\cdot G_{n-3}
\]
they get $\al_{\emph{SG}}\approx 0.38$. Our technique easily yields the exact values
\[\alpha_{\emph{Shu}}=\arctan\left(1/{\sqrt{135+78 \sqrt{3}}}\right)
\approx 0.0607719,\]
and
\[\alpha_{\emph{SG}}=\arctan\sqrt{\frac{1}{3}\left(2 \sqrt{3}-3\right)}
\approx 0.374734.\]
\end{remark}

\section{Optimal parabola inclusion in a family of multistep methods}\label{section8}

In the previous section we demonstrated how one can find the optimal sector in a family
of stability regions of multistep methods. Here we show that
the same algebraic approach allows us to replace the sector with more general shapes:
we use again the multistep family \eqref{modelexample} as a test example and determine the optimal stability region that contains the largest parabola. The motivation for considering 
the shape of a parabola
  comes from \cite{hundsdorferruuth} 
(``for advection-diffusion equations, stability within
a parabola\footnote{\label{footnote3}
Similarly to Footnote \ref{footnote2}, 
an analogous typo is present in \cite[Section 2.3]{hundsdorferruuth} when  
the notion of ``stability within 
a parabola'' is defined. There we should have again $\tan$ instead of $\arctan$, that is,
\[
\xi\le 0, \quad |\eta^2/\xi|\le \tan(\beta)\quad \text{ with some angle } 0< \beta\le\pi/2.
\]
} can be more relevant than for a wedge''),  or from \cite[Sections 3--4]{calvo} 
(where linearly implicit Runge--Kutta methods are developed for the numerical integration
of semidiscrete equations originating from spatial discretizations of 
PDEs of advection-reaction-diffusion type). 

With $P_{\beta_1,\beta_0}$ and $\Ss_{\beta_1,\beta_0}$ defined in 
\eqref{Pbetadef}--\eqref{Sbetadef}, we are now looking for the largest possible 
$m>0$ such that 
the stability region of a suitable member of the family \eqref{modelexample} 
contains the parabola 
\begin{equation}\label{paraboladefinition}
{\mathcal{P}_{m}} 
:=\{ (\xi, \eta) \in\mathbb{R}^2: \xi\le 0, \,  \eta\in\mathbb{R}, \, \eta^2\le m|\xi|\},
\end{equation}
that is, the inclusion
\begin{equation}\label{parabolainclusion}
{\mathcal{P}_{m}} \subset \Ss_{\beta_1,\beta_0}
\end{equation}
holds. Clearly, we need
$\overset{\circ}{\text{A}}$-stability again to have \eqref{parabolainclusion} with some $m>0$, so from now on, by Lemma \ref{Wlemma}, we can assume  that
$(\beta_1,\beta_0) \in W$ (see Figure \ref{fig:wedge}). 

In Appendix \ref{subsection81}
we apply a simple geometric argument: we first formulate the RLCs for the members of the multistep family as implicit curves $\{(\xi, \eta)\in \mathbb{R}^2 : F_{\beta_1, \beta_0}(\xi,\eta)=0\}$, then invoke the 
notion of discriminant \cite{gelfand} to construct a polynomial in $m$
(and depending on the parameters
$\beta_1$ and $\beta_0$) whose suitable root can yield the optimal value
$\widetilde{m}_{\text{opt}}$ in \eqref{parabolainclusion}. 
The simple observation is the same as the one used in Section \ref{section23} (or in Section \ref{diskinclusionsection}): the 
optimal inscribed object (now a parabola) touches the boundary of the optimal stability region.

Based on this technique and by using \textit{Mathematica}, we conjecture that 
the parameter values $\beta_1=1/5$ and $\beta_0=37/40$ give  
 $\widetilde{m}_{\text{opt}}=6/5$. In Appendix \ref{subsection82} we use a uniqueness argument to 
rigorously prove this conjecture. We emphasize that, similarly to Appendix \ref{optimalmethodW2},  no RLCs are
involved in this uniqueness proof; the RLCs are used only as auxiliary objects
to conjecture the optimum. Given the complexity of intermediate calculations, 
it is again surprising that the final result  $\widetilde{m}_{\text{opt}}$ is a simple rational number. In summary, we have the following theorem.

\begin{theorem}\label{parabolathm}
Suppose that $(\beta_1,\beta_0)\in W$. Then the largest $m>0$ such that
\eqref{parabolainclusion} holds is $m\equiv \widetilde{m}_{\text{opt}}:=6/5$.
\end{theorem}
\begin{remark}
The authors of \cite{hundsdorferruuth} observe that ``for the methods considered in this paper, a large angle $\alpha$ will correspond to a large $\beta$'' (with $\alpha$ and $\beta$ interpreted in our Footnotes \ref{footnote2} and \ref{footnote3}). According to our results, the optimal $(\beta_1,\beta_0)$ parameter pairs $(3/8, 3/4)$ 
and $(1/5, 37/40)$---determining the stability regions with the largest inscribed sector and parabola, respectively---do not coincide, although they are both located on the right boundary of $W$ in Figure \ref{fig:wedge} (see also Remark \ref{remark12}).
\end{remark}




\appendix 

\section{Appendix}\label{appendixA}

\subsection{The proof of Lemma \ref{Wlemma}}\label{sectionsomeneccond}

\begin{proof}
Let us fix some $(\beta_1, \beta_0)\in\mathbb{R}^2$. For $\xi<0$, $P_{\beta_1,\beta_0}(\zeta,\xi,0)=0$ is equivalent to $\text{LHS}(\zeta)=\text{RHS}(\zeta)$ with
\[
\text{LHS}_{\beta_1,\beta_0}(\zeta):=\beta _0 \zeta ^3+\beta _1 \zeta ^2
-\left(3 \beta _0+2 \beta _1-3\right) \zeta +2 \beta
   _0+\beta _1-\frac{3}{2}
\]
and
$
\text{RHS}_\xi(\zeta):={(\zeta ^3-{3 \zeta ^2}/{4}-{1}/{4})}/{\xi }$.
Clearly, if $|\xi|$ is large enough, the coefficients of the $\text{RHS}$ polynomial can be arbitrarily close to $0$. So by 
 the fact that the roots of a polynomial are continuous functions of its coefficients, we get that 
``the 
$\zeta_j$ roots of $\text{LHS}(\zeta)=\text{RHS}(\zeta)$  can be made arbitrarily close to those of  
$\text{LHS}(\zeta)=0$ by choosing $|\xi|$ large''. To make the previous ``statement''
precise, we distinguish two cases according to whether the leading coefficient of 
$\text{LHS}$ 
vanishes or not: 
for $\beta_0=0$, the $\text{LHS}$ polynomial has at most
two roots, whereas the difference $\text{LHS}-\text{RHS}$ has three.

Case I: $\beta_0\ne 0$. By the above statement we easily see that if 
$\text{LHS}_{\beta_1,\beta_0}(\cdot)\notin \vn$, then 
$P_{\beta_1,\beta_0}(\cdot,\xi,0)\notin \svn$ for $|\xi|$ large enough. We now show that
\begin{equation}\label{notinW2implies}
(\beta_1,\beta_0)\notin W \implies \text{LHS}_{\beta_1,\beta_0}(\cdot)\notin \vn.
\end{equation}
So let us suppose in the rest of Case I that $(\beta_1,\beta_0)\notin W$ and $\beta_0\ne 0$.

Case I$_\text{a}$. First we check the case  when 
$\cc \text{LHS}_{\beta_1,\beta_0}(\cdot)=0$. Then 
\[
\text{LHS}_{\beta_1,\beta_0}(\zeta)=\zeta/4 \left[\left(2 \beta _1-3\right) \zeta ^2-4 \beta _1 \zeta +2 \beta _1-3\right],
\]
and, since now $2 \beta _1-3\ne 0$, we can apply Theorem 
\ref{vnthm} to the above polynomial in $[ \cdots ]$: due to $[ \cdots ]^\mathbf{r}\equiv 0$ we have that
$[ \cdots ]\in\vn$ if and only if $\zeta\mapsto [ \cdots ]'=2 \left(2 \beta _1-3\right) \zeta -4 \beta _1\in\vn$. But we directly see that this last linear polynomial $\notin \vn$, because 
$(\beta_1,\beta_0)\notin W$ and $\cc \text{LHS}_{\beta_1,\beta_0}(\cdot)=0$
imply $\beta_1>3/4$.

Case I$_\text{b}$. The conditions $\cc \text{LHS}_{\beta_1,\beta_0}(\cdot)\ne 0 \ne \beta_0$ mean
that we can apply Theorem \ref{vnthm} to $\text{LHS}_{\beta_1,\beta_0}(\cdot)$. It is easy to
verify that $\left(\text{LHS}_{\beta_1,\beta_0}(\cdot)\right)^\mathbf{r}$ does not vanish identically, so 
$
\text{LHS}_{\beta_1,\beta_0}(\cdot)\in\vn 
$ if and only if 
\begin{equation}\label{neveroccurs}
\left|\lc \text{LHS}_{\beta_1,\beta_0}(\cdot)\right| > \left|\cc \text{LHS}_{\beta_1,\beta_0}(\cdot)\right| 
\text{ and } \left(\text{LHS}_{\beta_1,\beta_0}(\cdot)\right)^\mathbf{r} \in\vn.
\end{equation}
We show in Cases I$_\text{b1}$ and I$_\text{b2}$ below that \eqref{neveroccurs} never occurs. First
 we observe that the inequality constraint in \eqref{neveroccurs}
yields that $\lc \left(\text{LHS}_{\beta_1,\beta_0}(\cdot)\right)^\mathbf{r} \ne 0$.

Case I$_\text{b1}$. If $\cc \left(\text{LHS}_{\beta_1,\beta_0}(\cdot)\right)^\mathbf{r} = 0$, then 
the polynomial $\left(\text{LHS}_{\beta_1,\beta_0}(\zeta)\right)^\mathbf{r}$
has exactly two roots: $\zeta_1=0$ and
\[
\zeta_2=2-\frac{3}{2 \left(6 \beta _0+2 \beta _1-3\right)}+\frac{3}{2 \left(2 \beta _0+2 \beta
   _1-3\right)}.
\]
One directly checks that $\cc \left(\text{LHS}_{\beta_1,\beta_0}(\cdot)\right)^\mathbf{r} = 0$
and $(\beta_1,\beta_0)\notin W$ imply $|\zeta_2|>1$.

Case I$_\text{b2}$. If $\cc \left(\text{LHS}_{\beta_1,\beta_0}(\cdot)\right)^\mathbf{r} \ne 0$,
we apply Theorem \ref{vnthm} to get that the quadratic polynomial 
$\left(\text{LHS}_{\beta_1,\beta_0}(\cdot)\right)^\mathbf{r}\in\vn$ if and only if
either Case I$_{\text{b2}\alpha}$ or I$_{\text{b2}\beta}$ below occurs.

Case I$_{\text{b2}\alpha}$: when
$\left(\text{LHS}_{\beta_1,\beta_0}(\cdot)\right)^\mathbf{rr}\equiv 0$
and $\left[\left(\text{LHS}_{\beta_1,\beta_0}(\cdot)\right)^\mathbf{r}\right]'\in \vn$. 
In this case, however, the unique root of the polynomial $\left[\left(\text{LHS}_{\beta_1,\beta_0}(\cdot)\right)^\mathbf{r}\right]'$,
\[
\zeta_1=1-\frac{3}{4 \left(6 \beta _0+2 \beta _1-3\right)}+\frac{3}{4 \left(2 \beta _0+2 \beta
   _1-3\right)},
\]
has absolute value $>1$.

Case I$_{\text{b2}\beta}$: when 
$\left|\lc \left(\text{LHS}_{\beta_1,\beta_0}(\cdot)\right)^\mathbf{r}\right|>
\left|\cc \left(\text{LHS}_{\beta_1,\beta_0}(\cdot)\right)^\mathbf{r}\right|$ 
and
$\left(\text{LHS}_{\beta_1,\beta_0}(\cdot)\right)^\mathbf{rr}$ $\in \vn$. But then the
unique root of $\left(\text{LHS}_{\beta_1,\beta_0}(\cdot)\right)^\mathbf{rr}$ is
\[
\zeta_1=1-\frac{3 \left(2 \beta _0+2 \beta _1-3\right)}{24 \beta _0^2+32 \beta _1 \beta _0-36
   \beta _0+8 \beta _1^2-18 \beta _1+9},
\]
for which we again have $|\zeta_1|>1$, completing Case I.

Case II: $\beta_0=0$. Then
\[
P_{\beta_1, 0}(\zeta,\xi,0)=
\zeta ^3-\zeta ^2 \left(\beta _1 \xi +\frac{3}{4}\right)-\left(3-2 \beta _1\right) \xi \zeta -
\left(\beta _1 \xi -\frac{3 \xi }{2}+\frac{1}{4}\right),
\]
and the leading coefficient of this cubic polynomial is 1. For each fixed $\beta_1\in
\mathbb{R}$ we see that at least one of its coefficients is unbounded as $\xi\to-\infty$,
so (by Vieta's formulae) at least one of its roots $\zeta(\xi)$ is unbounded as $\xi\to-\infty$. Hence 
$(-\infty,0)\times \{0\} \subset \Ss_{\beta_1, 0}$ cannot hold.
\end{proof}





\subsection{The proof of Theorem \ref{wedgethm}}\label{optimalmethodW2}

\begin{proof}
In the proof we suppose $m>0$ and, due to Lemma \ref{Wlemma}, that
 $(\beta_1,\beta_0)\in W$.\\

\noindent \textbf{Step 1.} Let us apply the same ideas as in Section \ref{sectionsomeneccond} but along the ray $\eta=-m\xi$. 
For $\xi<0$ we consider the roots of $P_{\beta_1,\beta_0}(\cdot,\xi, -m\xi)$, and get 
that
\[
\text{MLHS}_{\beta_1,\beta_0, m}(\cdot)\notin \vn \implies 
P_{\beta_1,\beta_0}(\cdot,\xi, -m\xi)\notin \svn
\]
for some $|\xi|$  large enough,  
where the corresponding ``modified left-hand side'' is defined as
\[
\text{MLHS}_{\beta_1,\beta_0, m}(\zeta):=\text{LHS}_{\beta_1,\beta_0}(\zeta)-\frac{3}{2} i m \zeta ^2,
\]
and we have also taken into account that $\lc \text{MLHS}_{\beta_1,\beta_0, m}(\cdot)=\beta_0 \ne 0$.
(The corresponding ``modified right-hand side'' would be the same $\text{RHS}_\xi(\zeta)$ 
as in Section \ref{sectionsomeneccond}.) 
Hence if the inclusion \eqref{inclusion} holds with some $m>0$, then
$\text{MLHS}_{\beta_1,\beta_0, m}(\cdot)\in \vn$.\\

\noindent \textbf{Step 2.} In this step we derive a necessary condition for 
$\text{MLHS}_{\beta_1,\beta_0, m}(\cdot)\in \vn$. 
First, one simply checks via Theorem \ref{vnthm} that 
$\cc \text{MLHS}_{\beta_1,\beta_0, m}(\cdot) = 0$,  
$\text{MLHS}_{\beta_1,\beta_0, m}(\cdot)\in \vn$ and $m>0$ cannot be simultaneously true. 
So we can suppose \[\lc \text{MLHS}_{\beta_1,\beta_0, m}(\cdot) \ne 0
\ne \cc \text{MLHS}_{\beta_1,\beta_0, m}(\cdot).\]  
We check that 
$\left(\text{MLHS}_{\beta_1,\beta_0, m}(\cdot)\right)^\mathbf{r}$ does not vanish
identically, and that 
\[\left|\lc \text{MLHS}_{\beta_1,\beta_0, m}(\cdot)\right| >
\left| \cc \text{MLHS}_{\beta_1,\beta_0, m}(\cdot)\right|.\] 
Then by Theorem \ref{vnthm} we have that
\[
\text{MLHS}_{\beta_1,\beta_0, m}(\cdot)\in \vn \Longleftrightarrow 
\left(\text{MLHS}_{\beta_1,\beta_0, m}(\cdot)\right)^\mathbf{r} \in \vn. 
\] 
Now we see that 
\[
\lc \left(\text{MLHS}_{\beta_1,\beta_0, m}(\cdot)\right)^\mathbf{r} \ne 0
\ne \cc \left(\text{MLHS}_{\beta_1,\beta_0, m}(\cdot)\right)^\mathbf{r},
\]
and $\left(\text{MLHS}_{\beta_1,\beta_0, m}(\cdot)\right)^\mathbf{rr}$ does not
vanish identically. Thus Theorem \ref{vnthm} yields that
\[\left(\text{MLHS}_{\beta_1,\beta_0, m}(\cdot)\right)^\mathbf{r} \in \vn\] if and
only if
\begin{equation}\label{tempignore}
\left |\lc \left(\text{MLHS}_{\beta_1,\beta_0, m}(\cdot)\right)^\mathbf{r}\right| >
\left| \cc \left(\text{MLHS}_{\beta_1,\beta_0, m}(\cdot)\right)^\mathbf{r}\right|
\end{equation}
and
\begin{equation}\label{eq27}
\left(\text{MLHS}_{\beta_1,\beta_0, m}(\cdot)\right)^\mathbf{rr} \in \vn.
\end{equation}
Clearly, $\deg \left(\text{MLHS}_{\beta_1,\beta_0, m}(\cdot)\right)^\mathbf{rr}\le 1$, and we directly confirm that \eqref{tempignore} implies that the degree is exactly 1.
From this we obtain that 
\eqref{tempignore} and \eqref{eq27} hold if and only if \eqref{tempignore} and
\[
\left| 1+i m+\frac{2 i m \beta _0}{4 \beta _0+2 \beta _1-3}-
\right.
\]
\begin{equation}\label{ineq28}
\left.
\frac{3(1+i m)
   \left[ \beta _0 \left(6 m^2+2\right)+\left(2 \beta _1-3\right)
   \left(m^2+1\right)\right]}{24 \beta _0^2+4 \beta _0 \left(8 \beta _1+3 m^2-9\right)+\left(2 \beta
   _1-3\right) \left(4 \beta _1+3 m^2-3\right)}\right| \leq 1
\end{equation}
hold. In particular, \eqref{tempignore} guarantees that the denominators appearing in
\eqref{ineq28} are non-zero, hence from now on we can restrict the parameters $(\beta_1,\beta_0)\in W$ to the set $(\beta_1,\beta_0)\in W\setminus L$ with
\begin{equation}\label{Ldef}
L:=\left\{(\beta_1,\beta_0) \in\mathbb{R}^2: \beta _0= \frac{3-2 \beta _1}{4}\right\},
\end{equation}
being the left edge of the wedge $W$; see Figure \ref{fig:wedgecornermethod}. 
\begin{figure}[H]
\subfigure{
\includegraphics[width=0.5\textwidth]{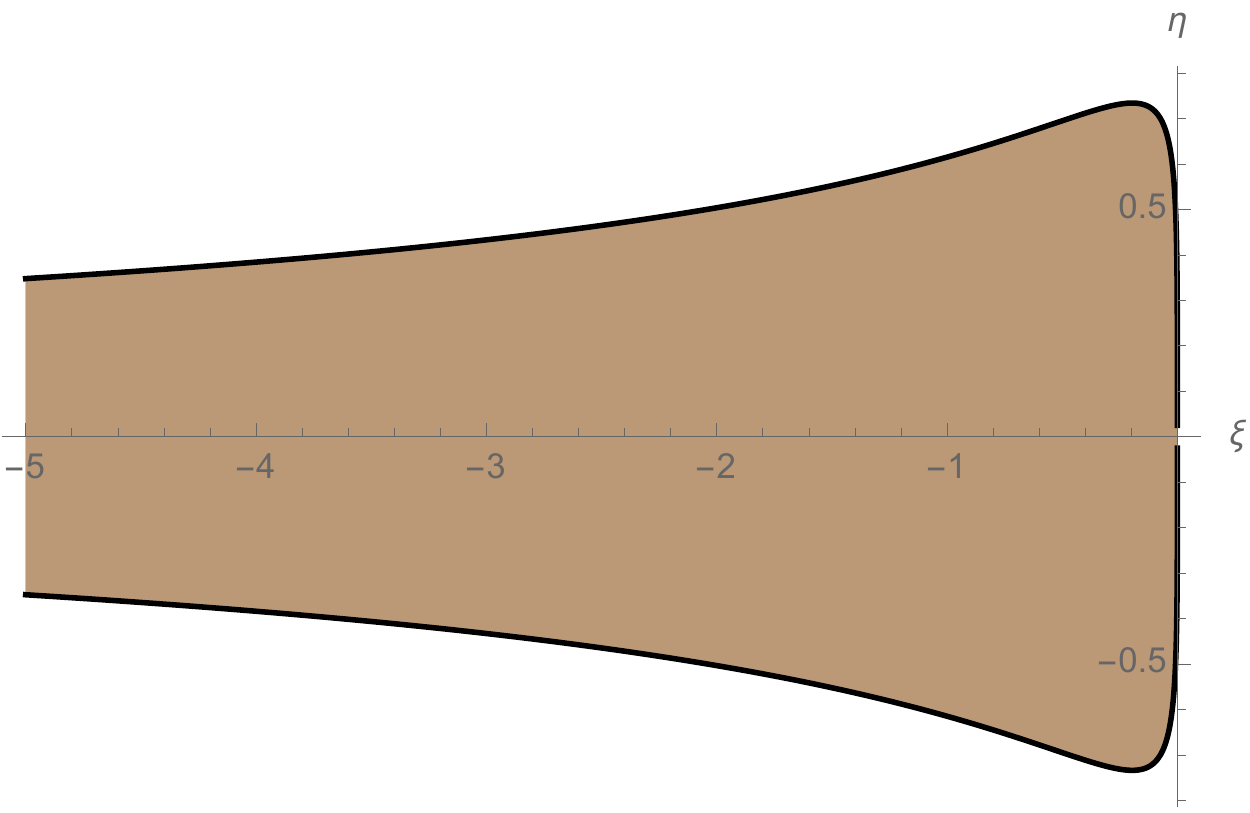}}
\subfigure{
\includegraphics[width=0.5\textwidth]{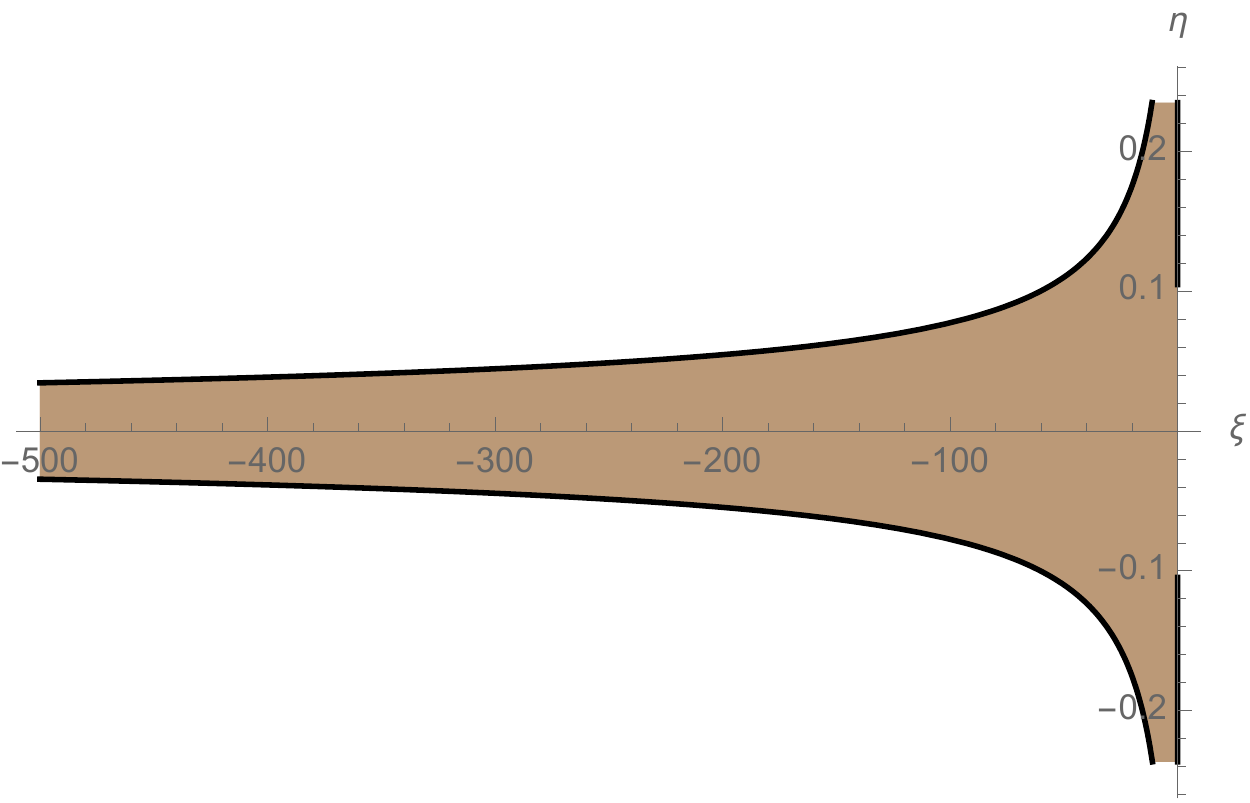}}
\caption{These figures show the stability region $\Ss_{3/4, 3/8}$ corresponding
to the method with $(\beta_1,\beta_0)=(3/4,3/8)$ (i.e., the vertex of the wedge in
Figure \ref{fig:wedge}). Such methods with $(\beta_1,\beta_0)\in
W\cap L$ (see \eqref{Ldef}) are $\overset{\circ}{\text{A}}$-stable, but 
${\mathcal{A}}_m\subset \Ss_{\beta_1,\beta_0}$ (see \eqref{inclusion}) does not hold with any $m>0$. \label{fig:wedgecornermethod}}
\end{figure}
By defining 
\[
C_4:= -9 \left(4 \beta _0+2 \beta _1-3\right)^2,
\]
\[
C_2:=
2 \big[864 \beta _0^4+864 \left(2 \beta _1-3\right) \beta _0^3+288 \left(4 \beta _1^2-13 \beta _1+10\right) \beta _0^2+
\]
\[4 \left(80 \beta _1^3-420 \beta _1^2+684 \beta
   _1-351\right) \beta _0+\left(3-2 \beta _1\right){}^2 \left(8 \beta _1^2-36 \beta _1+27\right)\big],
\]
\[
C_0:=-3 \left(4 \beta _0+2 \beta _1-3\right){}^2 \left(8 \beta _0+8 \beta _1-9\right)
\]
and
\[
Q_{\beta_1,\beta_0}(m):=C_4 m^4+C_2 m^2+C_0,
\]
it is easily verified after some factorization and simplification that 
\[
\eqref{tempignore} \text{ and } \eqref{ineq28} \Longleftrightarrow 
\eqref{tempignore} \text{ and } Q_{\beta_1,\beta_0}(m)\ge 0.
\]
In particular, $\text{MLHS}_{\beta_1,\beta_0, m}(\cdot)\in \vn$ implies 
$Q_{\beta_1,\beta_0}(m)\ge 0$.\\

\noindent \textbf{Step 3.} We see that $C_4<0$ and $C_0\ge 0$ for $(\beta_1,\beta_0)\in W\setminus L$, hence 
we can denote the largest real root of
the polynomial $Q_{\beta_1,\beta_0}(\cdot)$ by $m^*(\beta_1,\beta_0)\in [0,+\infty)$.
Consequently, if $\text{MLHS}_{\beta_1,\beta_0, m}(\cdot)\in \vn$, then $m\le m^*(\beta_1,\beta_0)$. 
We now conjecture (by using \textit{Mathematica}'s {\texttt{Maximize}} 
command, for example) that 
\begin{equation}\label{m*conjecture}
m^*(\beta_1,\beta_0)\le \frac{1}{2}\quad\text{for}\quad (\beta_1,\beta_0)\in W\setminus L,
\end{equation}
and $m^*(\beta_1,\beta_0)=1/2$ occurs precisely for
$(\beta _1,\beta_0)=({3}/{8}, {3}/{4})$, see Figure \ref{fig:mstar}.
With this conjectured optimal $m^*$ value, we can prove 
\eqref{m*conjecture} and the uniqueness property in an elementary way.

By introducing the shifted variable $M:=m-1/2$, we rewrite $Q_{\beta_1,\beta_0}(m)$ 
as 
\begin{equation}\label{Mpolydef}
\sum_{j=0}^4 \widehat{C}_j(\beta_1,\beta_0) \,M^j.
\end{equation}
 Then we check that
\[
(\beta_1,\beta_0)\in W\setminus L \implies \widehat{C}_j(\beta_1,\beta_0)<0
\text{ for } 1\le j\le 4.
\]
Moreover, we have 
\[
\widehat{C}_0(\beta_1,\beta_0)\equiv 6912 \beta _0^4+768 \left(18 \beta _1-35\right) \beta _0^3+48 \left(192 \beta _1^2-880
   \beta _1+813\right) \beta _0^2+
\]
\[
40 \left(64 \beta _1^3-528 \beta _1^2+1062 \beta
   _1-621\right) \beta _0+\left(3-2 \beta _1\right)^2 \left(64 \beta _1^2-672 \beta
   _1+639\right),
\]

\[
(\beta_1,\beta_0)\in W\setminus L \implies \widehat{C}_0(\beta_1,\beta_0)\le 0
\]
and 
\[
\left[(\beta_1,\beta_0)\in W\setminus L \text{\  and \ } \widehat{C}_0(\beta_1,\beta_0)=0\right]
\Longleftrightarrow (\beta_1,\beta_0)=(3/8,3/4).
\]
On the one hand, these mean that \eqref{Mpolydef} is negative for $M>0$ and $(\beta_1,\beta_0)\in W\setminus L$. On the other hand, for $M=0$ the polynomial \eqref{Mpolydef} is zero
if and only if $(\beta_1,\beta_0)=(3/8,3/4)$.
\begin{figure}[H]
\centerline{\includegraphics[width=.6\textwidth]{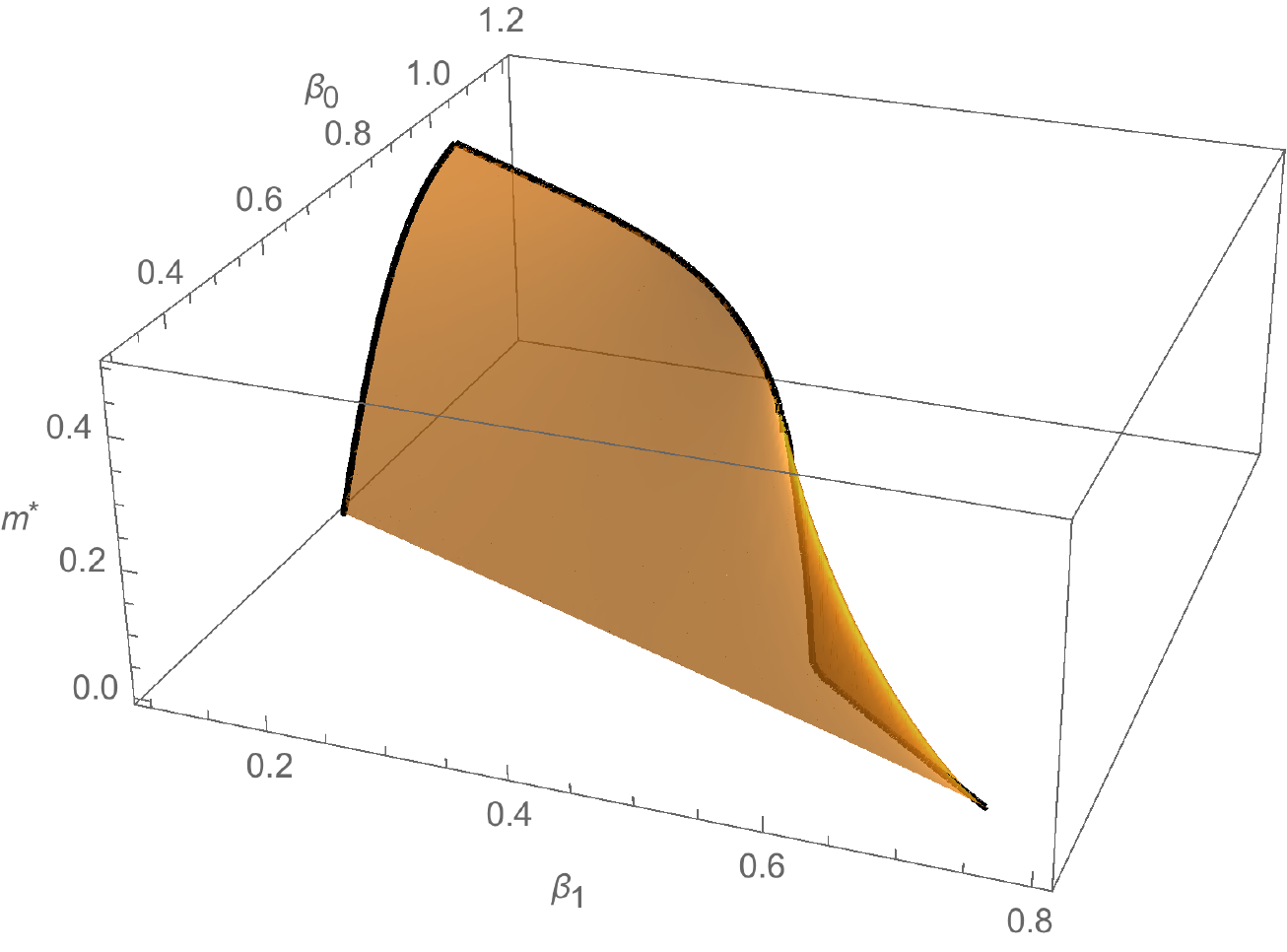}}
\caption{The function $m^*$ defined in Step 3 in Appendix \ref{optimalmethodW2}.
Its maximum value is located at $(\beta_1,\beta_0, m^*)=(3/8, 3/4,1/2)$.}
\label{fig:mstar}      
\end{figure}

Therefore we have proved that if \eqref{inclusion} holds with some $m>0$, then
$m\le 1/2$; and if $m=1/2$ is possible at all, then $(\beta _1,\beta_0)=({3}/{8}, {3}/{4})$.\\

\noindent \textbf{Step 4.} In this final step we show that $m=1/2$ in \eqref{inclusion} can be
achieved, by showing that ${\mathcal{A}}_{1/2}\subset \Ss_{{3}/{8}, {3}/{4}}$, that is,
\begin{equation}\label{A1/2implication}
(\xi,\eta)\in{\mathcal{A}}_{1/2}\implies P_{3/8,3/4}(\cdot,\xi, \eta) \in \svn.
\end{equation}
Let us fix such a pair $(\xi,\eta)$. One sees that 
\[
\left|\lc P_{3/8,3/4}(\cdot,\xi, \eta)\right|
 > \left|\cc P_{3/8,3/4}(\cdot,\xi, \eta)\right|,
\]
and in the $\cc P_{3/8,3/4}(\cdot,\xi, \eta) = 0$ case \eqref{A1/2implication} is easily 
verified to hold. 
Otherwise, if $\cc \ne 0$, we check that
$\left(P_{3/8,3/4}(\cdot,\xi, \eta)\right)^\mathbf{r}$ does not vanish identically, so by 
Theorem \ref{svnthm} we have that $P_{3/8,3/4}(\cdot,\xi, \eta)\in \svn$ if and only if
\begin{equation}\label{3/83/4firstreduced2}
\left(P_{3/8,3/4}(\cdot,\xi, \eta)\right)^\mathbf{r}\in \svn. 
\end{equation}
We have that
$\lc \left(P_{3/8,3/4}(\cdot,\xi, \eta)\right)^\mathbf{r}\ne 0$. Moreover,
$\cc \left(P_{3/8,3/4}(\cdot,\xi, \eta)\right)^\mathbf{r}= 0$ for $\xi=-2$ or 
$\xi=-2/3$, in which cases \eqref{3/83/4firstreduced2} holds.
So we can suppose from now on that $\cc \left(P_{3/8,3/4}(\cdot,\xi, \eta)\right)^\mathbf{r}\ne 0$.
Then one proves that 
\[
\lc \left(P_{3/8,3/4}(\cdot,\xi, \eta)\right)^\mathbf{rr}=
\]
\[
-\frac{81 \eta ^2 \xi ^2}{256}-\frac{27 \eta ^2 \xi }{64}-\frac{9 \eta ^2}{64}+\frac{81 \xi ^4}{512}-\frac{783 \xi ^3}{512}+\frac{441 \xi ^2}{128}-\frac{423 \xi
   }{128}+\frac{27}{32}\ne 0,
\]
so $\deg\left(P_{3/8,3/4}(\cdot,\xi, \eta)\right)^\mathbf{rr}=1$. Hence, by using
Theorem \ref{svnthm} again, we get that \eqref{3/83/4firstreduced2} holds if and only if
\[
\left|\lc \left(P_{3/8,3/4}(\cdot,\xi, \eta)\right)^\mathbf{r}\right|  > 
\left|\cc \left(P_{3/8,3/4}(\cdot,\xi, \eta)\right)^\mathbf{r}\right|
\text{ and } 
\left(P_{3/8,3/4}(\cdot,\xi, \eta)\right)^\mathbf{rr}\in \svn
\]
hold. Finally, we check that these last two conditions are satisfied for any $(\xi,\eta)\in{\mathcal{A}}_{1/2}$ pair not excluded earlier during the case separations.
\end{proof}

\begin{figure}[H]
\subfigure{
\includegraphics[width=0.5\textwidth]{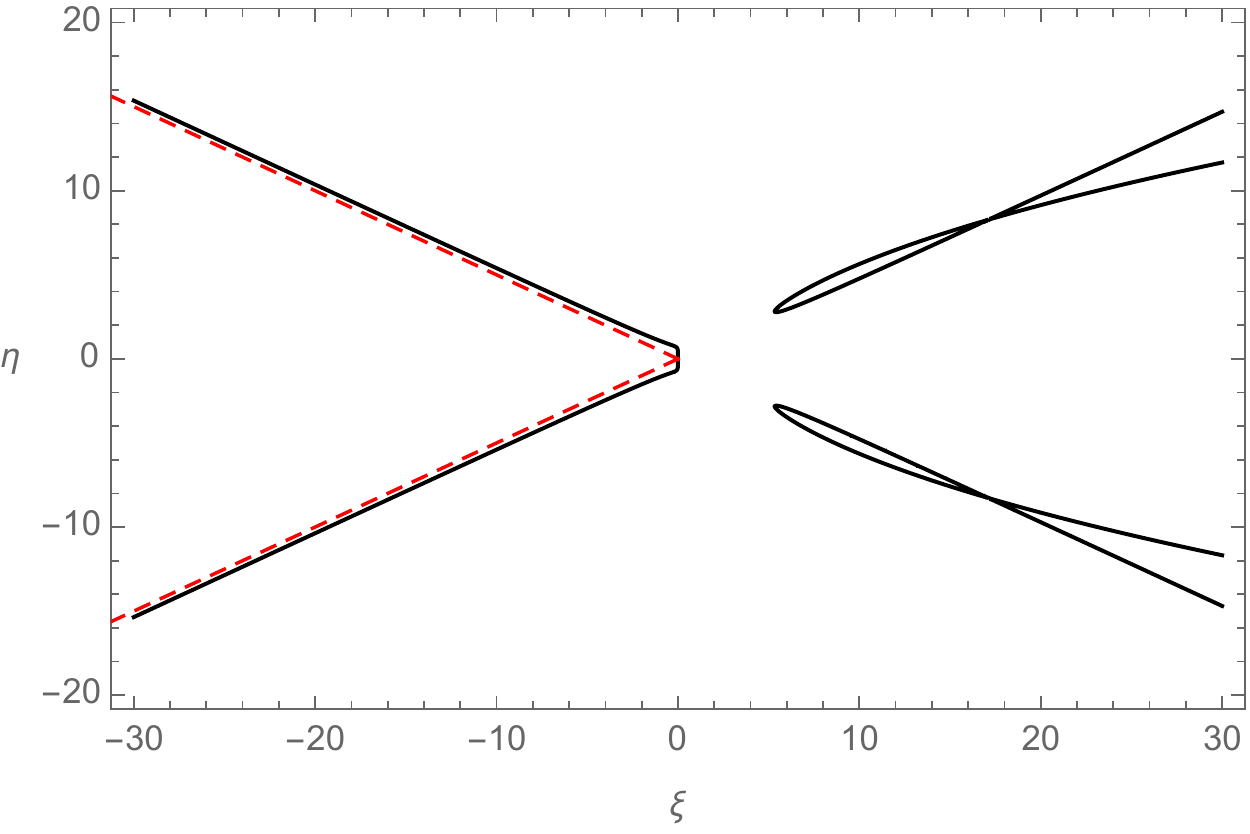}}
\subfigure{
\includegraphics[width=0.5\textwidth]{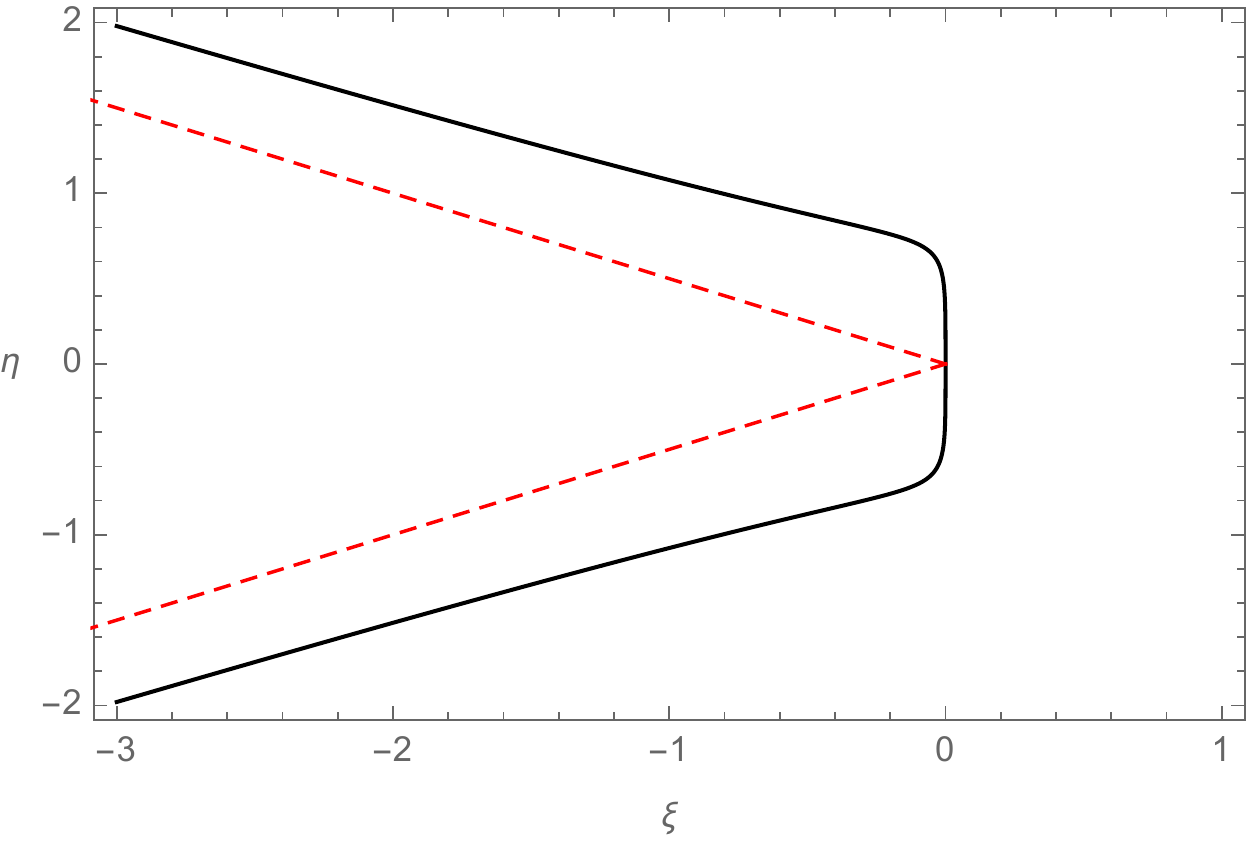}}
\caption{The implicit curve $\{(\xi, \eta)\in \mathbb{R}^2 : F_{\text{opt}}(\xi,\eta)=0\}$ 
(see Remark \ref{stabregioncoincideswithRLC}), being the boundary of the optimal
stability region $\Ss_{3/8, 3/4}$ in the left half-plane, 
 is shown in the left figure and
a close-up in the right figure in black. 
The dashed red lines represent the boundary of the largest infinite sector ${\mathcal{A}}_{1/2}$ that can be included
in the stability region.\label{fig:optimalIMEXmethod}}
\end{figure}

\begin{remark}\label{stabregioncoincideswithRLC}
By defining 
\[
F_{\text{opt}}(\xi,\eta):=12 \eta ^4 (3 \xi +2)^2-
\]
\[
3 \eta ^2 \xi  \left(9 \xi ^3+192 \xi ^2-620 \xi +368\right)+16 \xi  \left(3 \xi ^2-7 \xi +6\right)^2
\]
and applying Theorem \ref{svnthm}, it is straightforward to show  
(cf.~Step 4 in the above proof) that  
\[
\Ss_{3/8, 3/4}=\{(\xi, \eta) \in \mathbb{R}^2 : \xi\le 0,\, \eta\in\mathbb{R}, \, F_{\text{opt}}(\xi,\eta)\le 0\},
\]
see Figure \ref{fig:optimalIMEXmethod} and cf.~Remark \ref{remark12}.
\end{remark}

\section{Appendix}\label{appendixB}

\subsection{Locating the candidate optimum for Theorem \ref{parabolathm}}\label{subsection81}

For a given $(\beta_1, \beta_0)\in W$ pair, we can represent the RLC of 
the corresponding multistep method of the family \eqref{modelexample} as an implicit curve of the form
\begin{equation}\label{RLCfamilyimplicitform}
\{(\xi, \eta)\in \mathbb{R}^2 : \xi\le 0,\, \eta\in\mathbb{R},\,F_{\beta_1, \beta_0}(\xi,\eta)=0\}
\end{equation}
 by using the transformations in Section \ref{section23} as follows. First we perform the substitution
$\zeta\mapsto \frac{i-t}{i+t}$  in the polynomial \eqref{Pbetadef}, then
eliminate $t\in\mathbb{R}$ by taking the resultant of the real and imaginary parts of 
$P_{\beta_1,\beta_0}\left(\frac{i-t}{i+t},\xi,\eta\right)$. The resulting polynomial
can be factored to get
$
2^{34}\cdot 9 \cdot \left(1-\beta _0 \xi \right)^6\cdot  F_{\beta_1, \beta_0}(\xi,\eta)$;
the normalization with $F_{\beta_1, \beta_0}(0, 1)=9$ has been used
to make this polynomial $F_{\beta_1, \beta_0}$  unique. 
The term $\left(1-\beta _0 \xi \right)^6$
(cf.~the leading coefficient of $P_{\beta_1,\beta_0}(\cdot,\xi,\eta)$) does not vanish now
due to $\xi\le 0$ and
$(\beta_1, \beta_0)\in W$, hence \eqref{RLCfamilyimplicitform} is obtained. 
We are not going to display the polynomial $F_{\beta_1, \beta_0}(\xi,\eta)$ 
explicitly:  it contains 82 terms in its expanded
form and its degree in the variables/parameters 
$(\xi,\eta, \beta_1, \beta_0)$ is $(6, 4, 4 , 4)$.

Now supposing that the RLC  
\eqref{RLCfamilyimplicitform} describes the boundary of the stability region
of the multistep method determined by the given  pair
$(\beta_1, \beta_0)$, 
it is reasonable to expect that, say, the upper branch of the
largest parabola inscribed in $\Ss_{\beta_1,\beta_0}$,
$
\{(\xi, \eta)\in \mathbb{R}^2 : \xi< 0,\, \eta>0,\,\eta^2=-m \xi\}$,
 touches the RLC \eqref{RLCfamilyimplicitform} at some finite point. 
In this case, the polynomial
\[
(-\infty,0)\ni \xi \mapsto F_{\beta_1, \beta_0}\left(\xi,\sqrt{-m \xi}\right)
\]
 has a multiple root there---it is indeed a polynomial, because in our situation 
$F_{\beta_1, \beta_0}(\xi,\eta)$ contains only even powers of $\eta$ (namely,
$\eta^2$ and $\eta^4$). Moreover, we now  have
\[
F_{\beta_1, \beta_0}\left(\xi,\sqrt{-m \xi}\right)=\xi\cdot \widetilde{Q}_{\beta_1,\beta_0, m}(\xi),
\]
where $\widetilde{Q}_{\beta_1,\beta_0, m}(\cdot)$ is a quartic polynomial. 
The existence of a multiple root of
$\widetilde{Q}_{\beta_1,\beta_0, m}(\cdot)$ implies that the discriminant of 
this polynomial (with respect to $\xi$), denoted by
$\widetilde{\Delta}_{\beta_1,\beta_0}(m)$, vanishes. \textit{Mathematica} yields that
\[
\widetilde{\Delta}_{\beta_1,\beta_0}(m)=-2^{13}\cdot 3^6\cdot m^2\left(9 \beta _0+4 \beta _1-6\right)^2  \times 
\]
\[
\left(64 \beta _1^4
   m^3+\ldots-4410\right)^2 \left(590976\, \beta _0^2 \beta _1^3 m^5+\ldots-24402696417\right),
\]
where the  ``$\ldots$'' symbols contain 57 and 228 terms, respectively. We see that 
the factor $9 \beta _0+4 \beta _1-6$ above is always positive in $W$. 
In this way we can determine the parameter $m$  of the
 largest parabola within the 
region bounded by the RLC for any fixed $(\beta_1,\beta_0)\in W$. 

\begin{remark}\label{remark12}
By setting $(\beta_1,\beta_0)=(3/8, 3/4)$ for example
(corresponding to the ``sector-optimal'' method in Section \ref{optimalmethodW2}),
we have that 
the  RLC in \eqref{RLCfamilyimplicitform}
is identical to
$3/16\cdot  F_{\text{opt}}(\xi,\eta)$ in Remark 
\ref{stabregioncoincideswithRLC},
implying that the RLC in the  left half-plane $\xi\le 0$ represents the boundary of the
stability region $\Ss_{3/8, 3/4}$. Now $\widetilde{\Delta}_{3/8, 
3/4}(m)$ can be written as
\[
-\frac{3^{23}}{2^{13}} (3 m-1)^2 m^2 (m+4)^4 (3 m+16)^2 \left(36 m^3+1362 m^2+343
   m-2116\right),
\]
from which we can prove that the largest parabola ${\mathcal{P}_{m}}$  contained in 
$\Ss_{3/8, 3/4}$ has $m\approx 1.11226$
(being the unique positive root of the polynomial $\{36, 1362, 343, -2116\}$).
\end{remark}

By 
studying the positive roots of $\widetilde{\Delta}_{\beta_1,\beta_0}(\cdot)$
as $(\beta_1,\beta_0)$ is varied within $W$, we can \textit{conjecture}
that the value of $m$ in \eqref{parabolainclusion} cannot be greater than $6/5$
for the family \eqref{modelexample}. Moreover, $m=6/5$ occurs only
for $\beta_1=1/5$ and $\beta_0=37/40$, and in this case the
RLC and the upper parabola branch touch each other at $(\xi,\eta)=(-10/7, 2\sqrt{3/7})$. Since the polynomial 
$\widetilde{\Delta}_{\beta_1,\beta_0}(m)$ is 
much more complicated than the corresponding 
polynomial $Q_{\beta_1,\beta_0}(m)$ in Appendix \ref{optimalmethodW2}, 
this time \textit{Mathematica} could not confirm in a reasonable amount of computing time
that the value $m=6/5$ is indeed the optimal one.

\subsection{The proof of optimality in Theorem \ref{parabolathm}}\label{subsection82}

However, once the unique optimum has been  conjectured properly, the proof 
of optimality 
becomes straightforward to complete.\\

\noindent \textbf{Step 1.} By assuming $(\beta_1, \beta_0)\in W$ throughout the step,
we show that the point $(\xi_0,\eta_0):=(-10/7, 2\sqrt{3/7})$ 
belongs to precisely one stability region in the family, by verifying that 
\[
P_{\beta_1,\beta_0}\left( \cdot , \xi_0, \eta_0\right)\in \svn 
\Longleftrightarrow (\beta_1, \beta_0)=
\left(1/5, 37/40\right).
\]
To see this, first we check that
$\lc P_{\beta_1,\beta_0}\left( \cdot, \xi_0, \eta_0\right)\ne 0$. 
Moreover, it is easily seen that $\cc P_{\beta_1,\beta_0}\left( \cdot, \xi_0, \eta_0\right)$
vanishes exactly for $\beta _1\le 23/40$ and $\beta _0=\left(67-40 \beta _1\right)/80$,
and in this case the polynomial $P_{\beta_1,\beta_0}\left( \cdot, \xi_0, \eta_0\right)/\xi$
has $\deg =2$ but $\notin \svn$, as a recursive application of Theorem \ref{svnthm} shows.
Then we can also prove that 
$\left(P_{\beta_1,\beta_0}\left( \cdot, \xi_0, \eta_0\right)\right)^\mathbf{r}$ does not vanish
identically, and that 
\[
\left|\lc P_{\beta_1,\beta_0}\left( \cdot, \xi_0, \eta_0\right)\right|>
\left|\cc P_{\beta_1,\beta_0}\left( \cdot, \xi_0, \eta_0\right)\right|.
\]
Thus, according to Theorem \ref{svnthm}, 
\[
P_{\beta_1,\beta_0}\left( \cdot , \xi_0, \eta_0\right)\in \svn 
\Longleftrightarrow \left(P_{\beta_1,\beta_0}\left( \cdot, \xi_0, \eta_0\right)\right)^\mathbf{r}
\in \svn.
\]
Now we repeat the above process 
with $\left(P_{\beta_1,\beta_0}\left( \cdot, \xi_0, \eta_0\right)\right)^\mathbf{r}$.
We prove that 
\[\left|\lc \left(P_{\beta_1,\beta_0}\left( \cdot, \xi_0, \eta_0\right)\right)^\mathbf{r}\right|
>\left|\cc \left(P_{\beta_1,\beta_0}\left( \cdot, \xi_0, \eta_0\right)\right)^\mathbf{r}\right|>0\]
 and that 
$\left(P_{\beta_1,\beta_0}\left( \cdot, \xi_0, \eta_0\right)\right)^\mathbf{rr}$ does not vanish identically, so by Theorem \ref{svnthm} we have that
\[
\left(P_{\beta_1,\beta_0}\left( \cdot, \xi_0, \eta_0\right)\right)^\mathbf{r}\in \svn
\Longleftrightarrow \left(P_{\beta_1,\beta_0}\left( \cdot, \xi_0, \eta_0\right)\right)^\mathbf{rr}
\in \svn.
\]
But $\left(P_{\beta_1,\beta_0}\left( \cdot, \xi_0, \eta_0\right)\right)^\mathbf{rr}$ is a linear 
polynomial (it is easily checked that it cannot be a constant polynomial), so its unique 
(non-real complex) root can be directly expressed: one sees that the absolute value of this
root is $\le 1$ if and only if
\[
\left(8 \beta _0+8 \beta _1-19\right) \left(120 \beta _0+40 \beta _1-39\right)
\left(483840000 \beta _0^4+967680000 \beta _1 \beta _0^3-\right.
\]
\[
 1989440000 \beta
   _0^3+645120000 \beta _1^2 \beta _0^2-2967744000 \beta _1 \beta _0^2+2890070400 \beta
   _0^2+
\]
\[
179200000 \beta _1^3 \beta _0-1404096000 \beta _1^2 \beta _0+2856374400 \beta _1
   \beta _0-1693045320 \beta _0+
\]
\[
\left.
17920000 \beta _1^4-214336000 \beta _1^3+673766400 \beta
   _1^2-792582600 \beta _1+301631887\right)\le 0.
\]
The product of the first two factors is strictly negative in $W$, and a standard 
constrained optimization computation
shows that the third factor is $\ge 0$ in $W$ if and only if 
$(\beta_1, \beta_0)=\left(1/5, 37/40\right)$, completing Step 1.\\

\noindent \textbf{Step 2.} Since ${\mathcal{P}_{m_1}}\subseteq {\mathcal{P}_{m_2}}$ 
 is equivalent to $0<m_1\le m_2$ (see \eqref{paraboladefinition}), and now
$\left| \eta_0^2/\xi_0\right|=6/5$,
the uniqueness property in the previous step implies that $m\ge 6/5$ in 
\eqref{parabolainclusion} can hold only for  
$(\beta_1, \beta_0)=\left(1/5, 37/40\right)$. In this step we verify that
\eqref{parabolainclusion} indeed holds with $m=6/5$ and
$(\beta_1, \beta_0)=\left(1/5, 37/40\right)$, that is, we show that
$P_{1/5, 37/40}\left( \cdot, \xi, \eta\right)\in \svn$  for any
$(\xi, \eta) \in {\mathcal{P}_{6/5}}$.

Let us pick and fix an arbitrary point $(\xi, \eta) \in {\mathcal{P}_{6/5}}$.
Then we easily see that 
\[\left| \lc P_{1/5, 37/40}\left( \cdot, \xi, \eta\right)\right| > 
\left| \cc P_{1/5, 37/40}\left( \cdot, \xi, \eta\right)\right|,\] and this
$\cc =0$ if and only if $\xi=-5/11$; in this case Theorem \ref{svnthm} tells us that 
$\zeta\mapsto P_{1/5, 37/40}\left( \zeta, -5/11, \eta\right) = 
\zeta  (125 \zeta ^2-132 i \zeta  \eta -58 \zeta -7)/88\in \svn$.  So for $\xi\ne -5/11$, again
by Theorem \ref{svnthm} we get that
\[
P_{1/5, 37/40}\left( \cdot, \xi, \eta\right)\in \svn
\Longleftrightarrow \left(P_{1/5, 37/40}\left( \cdot, \xi, \eta\right)\right)^\mathbf{r}
\in \svn,
\]
provided that $\left(P_{1/5, 37/40}\left( \cdot, \xi, \eta\right)\right)^\mathbf{r}$ does 
not vanish identically. But this non-vanishing condition is true because
\[
\left|\lc \left(P_{1/5, 37/40}\left( \cdot, \xi, \eta\right)\right)^\mathbf{r}\right|
>\left|\cc \left(P_{1/5, 37/40}\left( \cdot, \xi, \eta\right)\right)^\mathbf{r}\right|>0.
\] 
Moreover, since
$\left|\lc \left(P_{1/5, 37/40}\left( \cdot, \xi, \eta\right)\right)^\mathbf{rr}\right|$ is
also positive, the above with Theorem \ref{svnthm}  imply that 
\[
P_{1/5, 37/40}\left( \cdot, \xi, \eta\right)\in \svn
\Longleftrightarrow \left(P_{1/5, 37/40}\left( \cdot, \xi, \eta\right)\right)^\mathbf{rr}
\in \svn.
\]
The positivity of $\left|\lc \left(\ldots\right)^\mathbf{rr}\right|$ yields that 
$\left(P_{1/5, 37/40}\left( \cdot, \xi, \eta\right)\right)^\mathbf{rr}$, a 
$\deg =1$ polynomial, has a unique root. The absolute value of this (real or complex)
root is $\le 1$ if and only if $(3 \xi -10) (59 \xi -30)\cdot \widetilde{F}_{\text{opt}}(\xi,\eta)\le 0$,
where
\[
\widetilde{F}_{\text{opt}}(\xi,\eta):=
720 \eta ^4 (11 \xi +5)^2-
\]
\[
\eta ^2 \xi  \left(19575 \xi ^3+485696 \xi ^2-1009140 \xi
   +464400\right)+
240 \xi  \left(22 \xi ^2-49 \xi +30\right)^2.
\]
Now $(3 \xi -10) (59 \xi -30)>0$, and one checks that 
$\widetilde{F}_{\text{opt}}(\xi,\eta)\le 0$ for $(\xi, \eta) \in {\mathcal{P}_{6/5}}$,
completing Step 2.\\

\noindent \textbf{Step 3.} To complete the optimality proof, we finally show  
that
\[
P_{1/5, 37/40}\left( \cdot, \xi_0, \eta_0+\varepsilon\right)\notin \svn 
\text{  for any  }  \varepsilon\in (0,1),
\]
that is, we cannot have $m>6/5$ in \eqref{parabolainclusion}. We repeat 
the same two-step reduction process as above and get that
$\varepsilon\in (0,1)$ guarantees that 
\[
P_{1/5, 37/40}\left( \cdot, \xi_0, \eta_0+\varepsilon\right)\in\svn 
\Longleftrightarrow 
\left(P_{1/5, 37/40}\left( \cdot, \xi_0, \eta_0+\varepsilon\right)\right)^\mathbf{rr}\in\svn.
\]
But this last $\left(\ldots\right)^\mathbf{rr}\in\svn$ condition is equivalent to
\[
\frac{1120 \varepsilon  \left(27783 \varepsilon ^3+31752 \sqrt{21} \varepsilon ^2+1649620 \varepsilon +833776 \sqrt{21}\right)}{\left(1323 \varepsilon ^2+756 \sqrt{21} \varepsilon
   -50840\right)^2}\le 0,
\]
so it cannot hold for any $\varepsilon\in (0,1)$.

\begin{figure}[H]
\centerline{\includegraphics[width=.6\textwidth]{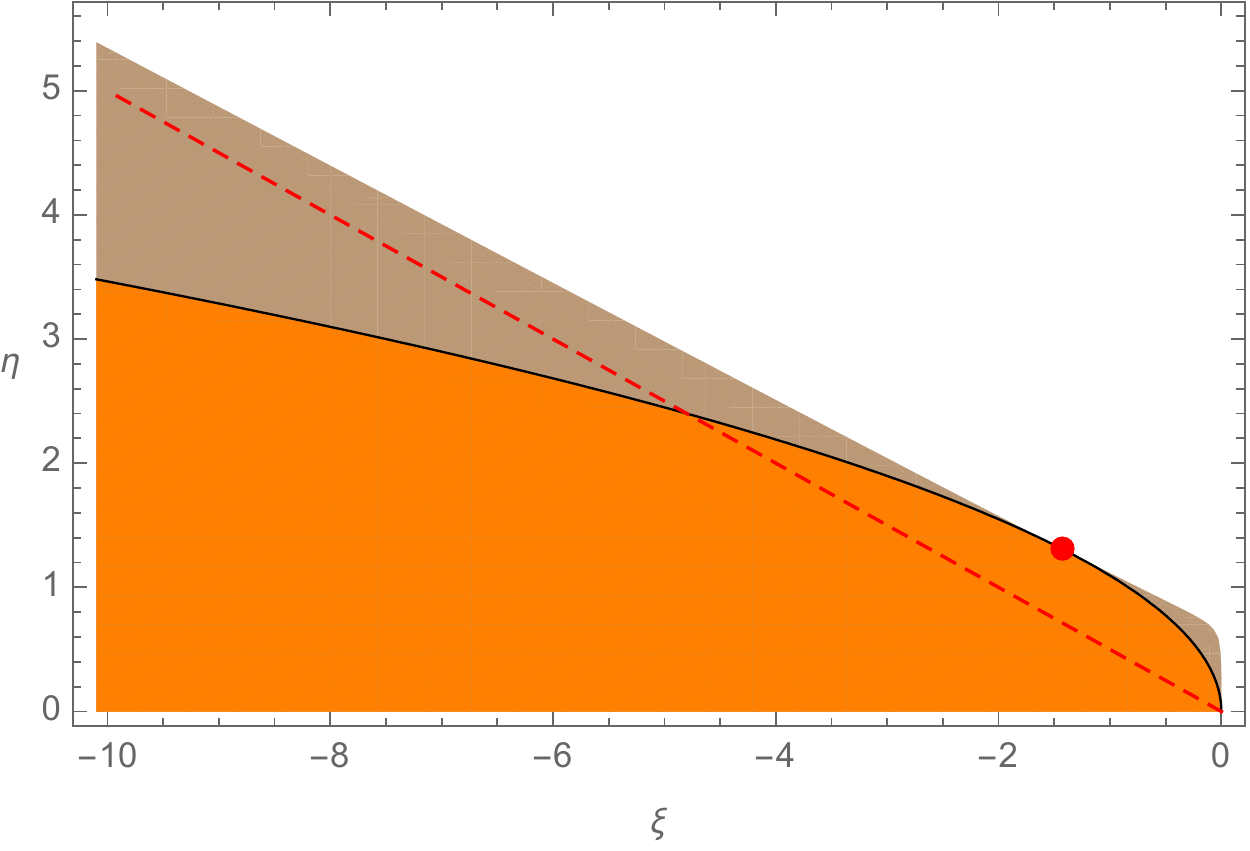}}
\caption{The figure shows the optimal stability region $\Ss_{1/5, 37/40}$ (brown)
within the family \eqref{modelexample} that contains the largest parabola
${\mathcal{P}_{6/5}}$ (orange), see Theorem \ref{parabolathm}. The point $(\xi_0,\eta_0)=(-10/7, 2\sqrt{3/7})$ 
is shown as a red dot. For comparison, the dashed red line from
Figure \ref{fig:optimalIMEXmethod} is also included here.}
\label{fig:otpimal_parabola}      
\end{figure}

\begin{remark}\label{closingremark16} 
In addition to the inequality $\widetilde{F}_{\text{opt}}(\xi,\eta)\le 0$  
in Step 2, we have that
$\widetilde{F}_{\text{opt}}(\xi,\eta)=0$ for $(\xi, \eta) \in {\mathcal{P}_{6/5}}$
if and only if $(\xi, \eta) =(\xi_0, \pm \eta_0)$. Moreover, 
$
\widetilde{F}_{\text{opt}}(\xi,\eta)\equiv 2000\cdot F_{1/5, 37/40}(\xi,\eta)
$
(see \eqref{RLCfamilyimplicitform}). On the other hand, by using the reduction process 
one can actually prove that 
\[
\{(\xi, \eta) \in \mathbb{R}^2 : \xi\le 0,\, \eta\in\mathbb{R}, \, \widetilde{F}_{\text{opt}}(\xi,\eta)\le 0\}=\Ss_{1/5, 37/40}.
\]  
These mean that the stability region boundary in the optimal case coincides with the 
corresponding RLC (in the left half-plane), 
and the boundary of the optimal inscribed parabola touches the 
stability region boundary in the open upper left half-plane at exactly one point, 
see Figure \ref{fig:otpimal_parabola} (and cf.~Remarks \ref{remarktouch} and \ref{remark12}).
\end{remark}

\end{document}